\documentclass[amstex,12pt, amssymb]{article}

\usepackage{mathtext}
\usepackage[cp1251]{inputenc}
\usepackage[T2A]{fontenc}
\usepackage[dvips]{graphicx}
\usepackage{amsmath}
\usepackage{amssymb}
\usepackage{amsxtra}
\usepackage{latexsym}
\usepackage{ifthen}

\newcounter{lemma}[section]

\newcounter{corollary}[section]

\newcounter{remark}[section]

\newcounter{theorem}[section]

\newcounter{proposition}[section]

\newcounter{example}

\numberwithin{equation}{section}

\begin{document}

\markboth{Z. KOVBA, E.~SEVOST'YANOV}{\centerline{ON MAPPINGS WITH
INVERSE INEQUALITIES ...}}

\def\cc{\setcounter{equation}{0}
\setcounter{figure}{0}\setcounter{table}{0}}

\overfullrule=0pt


\author{ZARINA KOVBA, EVGENY SEVOST'YANOV}

\title{\bf ON MAPPINGS WITH INVERSE
INEQUALITIES AND PRIME ENDS OF VARIABLE DOMAINS}

\date{\today}
\maketitle

\begin{abstract}
We consider mappings that satisfy inverse moduli inequalities of
Poletskii type, under which the images of the domain under the
mappings may change. It is proved that a classes of such mappings
are equicontinuous with respect to prime ends of some domain if the
majorant in the indicated modulus inequality is integrable.
\end{abstract}


\section{Introduction}

The present manuscript is devoted to the study of mappings
satisfying some estimates of the distortion of the modulus of
families of paths, the mapped domains under which may change, cf.
\cite{Cr}, \cite{MRSY}. The main subject of the study is the
equicontinuity of such families of mappings in terms of prime ends.
Some special cases of a similar problem have been studied by us
earlier, see e.g. \cite{ISS}, \cite{Sev$_1$}--\cite{Sev$_2$} and
\cite{SevSkv$_1$}. However, here we considered exclusively a fixed
domain under the mapping. Our goal is to remove this restriction and
consider families of mappings onto a sequence of domains. It should
be noted that, some related problems were previously studied by
G.~Suvorov and V.~Kruglikov (see, e.g.~\cite{Suv} and \cite{Kr}). In
particular, G.~Suvorov obtained some general results for plane
homeomorphisms of simply connected domains, including mappings with
a bounded Dirichlet integral (see~\cite{Suv}). In turn, V.~Kruglikov
slightly modified Suvorov's approach and thus obtained theorems on
boundary correspondence for mappings that are quasiconformal in the
mean (see~\cite{Kr}). Compared to Suvorov, some more general
theorems on the boundary correspondence of spatial quasiconformal
mappings were obtained. Unfortunately, the results from~\cite{Kr}
were not published in the wide press with full proofs. At the same
time, our goals at the present moment are not as deep as in the
mentioned studies, if we mean the construction of a whole geometric
theory of prime ends of sequences of domains. We will restrict
ourselves here to the case when the description of the behavior of
mappings is possible by using prime ends of the kernel of such a
sequence. Let us also note some classical and relatively modern
studies by other authors, see, e.g., \cite{Ad, ABBS, GU, KPR, KR,
Na$_2$}.

\medskip
Let us recall some definitions. A Borel function $\rho:{\Bbb
R}^n\,\rightarrow [0,\infty] $ is called {\it admissible} for the
family $\Gamma$ of paths $\gamma$ in ${\Bbb R}^n,$ if the relation
\begin{equation*}\label{eq1.4}
\int\limits_{\gamma}\rho (x)\, |dx|\geqslant 1
\end{equation*}
holds for all (locally rectifiable) paths $ \gamma \in \Gamma.$ In
this case, we write: $\rho \in {\rm adm} \,\Gamma .$ Let $p\geqslant
1,$ then {\it $p$-modulus} of $\Gamma $ is defined by the equality
\begin{equation*}\label{eq1.3gl0}
M_p(\Gamma)=\inf\limits_{\rho \in \,{\rm adm}\,\Gamma}
\int\limits_{{\Bbb R}^n} \rho^p (x)\,dm(x)\,.
\end{equation*}
We set $M(\Gamma):=M_n(\Gamma).$ Let $x_0\in {\Bbb R}^n,$
$0<r_1<r_2<\infty,$
\begin{equation}\label{eq1ED}
S(x_0,r) = \{ x\,\in\,{\Bbb R}^n : |x-x_0|=r\}\,, \quad B(x_0, r)=\{
x\,\in\,{\Bbb R}^n : |x-x_0|<r\}\end{equation}
and
\begin{equation}\label{eq1**}
A=A(x_0, r_1,r_2)=\left\{x\,\in\,{\Bbb R}^n:
r_1<|x-x_0|<r_2\right\}\,.\end{equation}
Given sets $E,$ $F\subset\overline{{\Bbb R}^n}$ and a domain
$D\subset {\Bbb R}^n$ we denote by $\Gamma(E,F,D)$ a family of all
paths $\gamma:[a,b]\rightarrow \overline{{\Bbb R}^n}$ such that
$\gamma(a)\in E,\gamma(b)\in\,F $ and $\gamma(t)\in D$ for $t \in
(a, b).$ Let $S_i=S(x_0, r_i),$ $i=1,2,$ where $S(x_0, r_i)$ are
defined in~(\ref{eq1ED}). If $f:D\rightarrow {\Bbb R}^n,$ $y_0\in f(
D)$ and $0<r_1<r_2<d_0=\sup\limits_{y\in f(D)}|y-y_0|,$ then by
$\Gamma_f(y_0, r_1, r_2)$ we denote the family of all paths $\gamma$
in $D$ such that $f(\gamma)\in \Gamma(S(y_0, r_1), S( y_0, r_2),
A(y_0,r_1,r_2)).$ Let $Q:{\Bbb R}^n\rightarrow [0, \infty]$ be a
Lebesgue measurable function. We say that {\it $f$ satisfies
Poletsky inverse inequality} at the point $y_0\in f(D),$ if the
relation
\begin{equation}\label{eq2*A}
M(\Gamma_f(y_0, r_1, r_2))\leqslant \int\limits_{A(y_0,r_1,r_2)\cap
f(D)} Q(y)\cdot \eta^n (|y-y_0|)\, dm(y)
\end{equation}
holds for any Lebesgue measurable function $\eta:
(r_1,r_2)\rightarrow [0,\infty ]$ such that
\begin{equation}\label{eqA2}
\int\limits_{r_1}^{r_2}\eta(r)\, dr\geqslant 1\,.
\end{equation}

\medskip
Recall that a mapping $f:D\rightarrow {\Bbb R}^n$ is called {\it
discrete} if the pre-image $\{f^{-1}\left(y\right)\}$ of each point
$y\,\in\,{\Bbb R}^n$ consists of isolated points, and {\it is open}
if the image of any open set $U\subset D$ is an open set in ${\Bbb
R}^n.$ A mapping $f:D\rightarrow {\Bbb R}^n$ is called {\it closed}
if $f(A)$ is closed in $f(D)$ whenever $A$ is closed in $D.$ Later,
in the extended space $\overline{{{\Bbb R}}^n}={{\Bbb
R}}^n\cup\{\infty\}$ we use the {\it spherical (chordal) metric}
$h(x,y)=|\pi(x)-\pi(y)|,$ where $\pi$ is a stereographic projection
$\overline{{{\Bbb R}}^n}$ onto the sphere
$S^n(\frac{1}{2}e_{n+1},\frac{1}{2})$ in ${{\Bbb R}}^{n+1},$ namely,
\begin{equation*}\label{eq3C}
h(x,\infty)=\frac{1}{\sqrt{1+{|x|}^2}}\,,\quad
h(x,y)=\frac{|x-y|}{\sqrt{1+{|x|}^2} \sqrt{1+{|y|}^2}}\,, \quad x\ne
\infty\ne y
\end{equation*}
(see \cite[Definition~12.1]{Va}). In what follows, ${\rm Int\,}A$
denotes the set of inner points of the set $A\subset \overline{{\Bbb
R}^n}.$ Recall that the set $U\subset\overline{{\Bbb R}^n}$ is a
neighborhood of the point $z_0,$ if $z_0\in {\rm Int\,}A.$

\medskip
The next definitions due to Caratheodory~\cite{Car}; cf.~\cite{GU,
KR} and~\cite{KrPa,Mikl} and~\cite{Suv}. Let $\omega$ be an open set
in ${\Bbb R}^k$, $k=1,\ldots,n-1$. A continuous mapping
$\sigma\colon\omega\rightarrow{\Bbb R}^n$ is called a {\it
$k$-dimensional surface} in ${\Bbb R}^n$. A {\it surface} is an
arbitrary $(n-1)$-dimensional surface $\sigma$ in ${\Bbb R}^n.$ A
surface $\sigma$ is called {\it a Jordan surface}, if
$\sigma(x)\ne\sigma(y)$ for $x\ne y$. In the following, we will use
$\sigma$ instead of $\sigma(\omega)\subset {\Bbb R}^n,$
$\overline{\sigma}$ instead of $\overline{\sigma(\omega)}$ and
$\partial\sigma$ instead of
$\overline{\sigma(\omega)}\setminus\sigma(\omega).$ A Jordan surface
$\sigma\colon\omega\rightarrow D$ is called a {\it cut} of $D$, if
$\sigma$ separates $D,$ that is $D\setminus \sigma$ has more than
one component, $\partial\sigma\cap D=\varnothing$ and
$\partial\sigma\cap\partial D\ne\varnothing$.

A sequence of cuts $\sigma_1,\sigma_2,\ldots,\sigma_m,\ldots$ in $D$
is called {\it a chain}, if:

(i) the set $\sigma_{m+1}$ is contained in exactly one component
$d_m$ of the set $D\setminus \sigma_m,$ wherein $\sigma_{m-1}\subset
D\setminus (\sigma_m\cup d_m)$; (ii)
$\bigcap\limits_{m=1}^{\infty}\,d_m=\varnothing.$

Two chains of cuts  $\{\sigma_m\}$ and $\{\sigma_k^{\,\prime}\}$ are
called {\it equivalent}, if for each $m=1,2,\ldots$ the domain $d_m$
contains all the domains $d_k^{\,\prime},$ except for a finite
number, and for each $k=1,2,\ldots$ the domain $d_k^{\,\prime}$ also
contains all domains $d_m,$ except for a finite number.

The {\it end} of the domain $D$ is the class of equivalent chains of
cuts in $D$. Let $K$ be the end of $D$ in ${\Bbb R}^n$, then the set
$I(K)=\bigcap\limits_{m=1}\limits^{\infty}\overline{d_m}$ is called
{\it the impression of} $K$. Following~\cite{Na$_2$}, we say that
the end $K$ is {\it a prime end}, if $K$ contains a chain of cuts
$\{\sigma_m\}$ such that
$$\lim\limits_{m\rightarrow\infty}M(\Gamma(C, \sigma_m, D))=0$$
for some continuum $C$ in $D.$ The following notation is used: the
set of prime ends corresponding to the domain $D,$ is denoted by
$E_D,$ and the completion of the domain $D$ by its prime ends is
denoted $\overline{D}_P.$

\medskip
Consider the following definition, which goes back to
N\"akki~\cite{Na$_2$}, cf.~\cite{KR}. The boundary of a domain $D$
in ${\Bbb R}^n$ is said to be {\it locally quasiconformal} if every
$x_0\in\partial D$ has a neighborhood $U$ that admits a
quasiconformal mapping $\varphi$ onto the unit ball ${\Bbb
B}^n\subset{\Bbb R}^n$ such that $\varphi(\partial D\cap U)$ is the
intersection of ${\Bbb B}^n$ and a coordinate hyperplane. The
sequence of cuts $\sigma_m,$ $m=1,2,\ldots ,$ is called {\it
regular,} if
$\overline{\sigma_m}\cap\overline{\sigma_{m+1}}=\varnothing$ for
$m\in {\Bbb N}$ and, in addition, $d(\sigma_{m})\rightarrow 0$ as
$m\rightarrow\infty.$ If the end $K$ contains at least one regular
chain, then $K$ will be called {\it regular}. We say that a bounded
domain $D$ in ${\Bbb R}^n$ is {\it regular}, if $D$ can be
quasiconformally mapped to a domain with a locally quasiconformal
boundary whose closure is a compact in ${\Bbb R}^n,$ and, besides
that, every prime end in $D$ is regular. Note that space
$\overline{D}_P=D\cup E_D$ is metric, which can be demonstrated as
follows. If $g:D_0\rightarrow D$ is a quasiconformal mapping of a
domain $D_0$ with a locally quasiconformal boundary onto some domain
$D,$ then for $x, y\in \overline{D}_P$ we put:
\begin{equation}\label{eq5M}
\rho(x, y):=|g^{\,-1}(x)-g^{\,-1}(y)|\,,
\end{equation}
where the element $g^{\,-1}(x),$ $x\in E_D,$ is to be understood as
some (single) boundary point of $D_0.$ The specified boundary point
is unique and well-defined by~\cite[Theorem~2.1, Remark~2.1]{IS},
cf.~\cite[Theorem~4.1]{Na$_2$}. It is easy to verify that~$\rho$
in~(\ref{eq5M}) is a metric on $\overline{D}_P.$ If $g_*$ is another
quasiconformal mapping of a domain $D_*$ with locally quasiconformal
boundary onto $D$, then the corresponding metric
$\rho_*(p_1,p_2)=|{\widetilde{g_*}}^{-1}(p_1)-{\widetilde{g_*}}^{-1}(p_2)|$
generates the same convergence and, consequently, the same topology
in $\overline {D}_P$ as $\rho_0$ because $g_0\circ g_*^{-1}$ is a
quasiconformal mapping of $D_*$ and $D_0$, which extends, by Theorem
4.1 in~\cite{Na$_2$}, to a homeomorphism between $\overline {D_*}$
and $\overline {D_0}$. In the sequel, this topology in $\overline
{D}_P$ will be called the {\it topology of prime ends}; the
continuity of mappings $F\colon \overline
{D}_P\rightarrow\overline{D^{\,\prime}}_P$ will be understood
relative to this topology.

\medskip
We say that the boundary of $D$ is {\it weakly flat} at a point
$x_0\in\partial D$ if, for every number $P > 0$ and every
neighborhood $U$ of the point $x_0,$ there is a neighborhood
$V\subset U$ such that $M_{\alpha}(\Gamma(E, F, D))\geqslant  P$ for
all continua $E$ and $F$ in $D$ intersecting $\partial U$ and
$\partial V.$ We say that the boundary $\partial D$ is weakly flat
if the corresponding property holds at every point of the boundary.

\medskip
Let $D_m,$ $m=1,2,\ldots ,$ be a sequence of domains in ${\Bbb
R}^n,$ containing a fixed point $A_0.$ If there exists a ball
$B(A_0, \rho),$ $\rho>0,$ belonging to all $D_m,$ then the {\it
kernel} of the sequence $D_m,$ $m=1,2,\ldots ,$ with respect to
$A_0$ is the largest domain $D_0$ containing $A_0$ and such that for
each compact set $E$ belonging to $D_0$ there is $N>0$ such that $E$
belongs to $D_m$ for all $m\geqslant N.$ A largest domain is one
which contains any other domain having the same property. A sequence
of domains $D_m,$ $m=1,2,\ldots ,$ converges to a kernel $D_0$ if
any subsequence of $D_m$ has $D_0$ as its kernel.

\medskip
Let $D_m,$ $m=1,2,\ldots ,$ be a sequence of domains which converges
to a kernel $D_0.$ Then $D_m$ will be called {\it regular} with
respect to $D_0,$ if $D_m\subset D_0$ for all $m\in {\Bbb N}$ and,
for every $P_0\in E_{D_0}$ there is a sequence $\sigma_k,$
$k=1,2,\ldots,$ with the following condition: if $d_k$ is a domain
in $P_0$ then there is $M=M(k)$ such that $d_k\cap D_m$ is a
non-empty connected set for every $m\geqslant M(k).$

\medskip
Given $E_1, E_2\subset\overline{{\Bbb R}^n}$ we set
$$h(E_1, E_2)=\inf\limits_{x\in E_1, y\in E_2}h(x, y)\,.$$
Given $\delta>0,$ domains $D, D_0\subset {\Bbb R}^n,$ $n\geqslant
2,$ a compact set $E$ in $D_0$ and Lebesgue measurable function
$Q:{\Bbb R}^n\rightarrow [0, \infty],$ $Q(y)\equiv 0$ outside $D_0,$
we denote by $\frak{F}^{E}_{Q, \delta}(D, D_0)$ a class of all open,
discrete and closed mappings $f:D\rightarrow D_0$ of a domain $D$
onto some domain $f(D),$ $E\subset f(D)\subset D_0,$ such that
(\ref{eq2*A})--(\ref{eqA2}) hold for every $y_0\in \overline{D_0}$
and, in addition, $h(f^{\,-1}(E),
\partial D)\geqslant \delta.$ The following statement was proved in
\cite{Sev$_2$} in a weaker form, when the image domain under the
mapping is fixed (cf.~\cite{Sev$_3$}, where this fact was
established on Riemannian manifolds under similar conditions). The
conditions on the function $Q,$ formulated in the indicated papers,
assume its integrability, which may also be weakened to the
integrability of this function over concentric spheres, see below.

\medskip
\begin{theorem}\label{th2}
{\it\, Let $\delta>0,$ let $D, D_0\subset {\Bbb R}^n$ be domains,
$n\geqslant 2,$ let $E$ be a compact set in $D_0$ and let $Q:{\Bbb
R}^n\rightarrow [0, \infty],$ $Q(y)\equiv 0$ outside $D_0,$ be a
Lebesgue measurable function. Assume that, 1) no connected component
of the boundary of the domain $D$ degenerates into a point, 2) $D$
has a weakly flat boundary and 3) $D_0$ is a regular domain. Let
$f_m\in \frak{F}^{E}_{Q, \delta}(D, D_0),$ $m=1,2,\ldots ,$ be a
sequence such that:

\medskip
4) every $f_m,$ $m=1,2,\ldots,$ has a continuous boundary extension
$f_m:\overline{D}\rightarrow \overline{D_0}_P,$

\medskip
5) the sequence of domains $f_m(D)$ is regular with respect
to~$D_0.$

\medskip
Assume that, for each point $y_0\in \overline{D_0}$ and for every
$0<r_1<r_2<r_0:=\sup\limits_{y\in D_0}|y-y_0|$ there is a set
$E_1\subset[r_1, r_2]$ of a positive linear Lebesgue measure such
that the function $Q$ is integrable with respect to
$\mathcal{H}^{n-1}$ over the spheres $S(y_0, r)$ for every $r\in
E_1.$

\medskip
Then the family $f_m,$ $m=1,2,\ldots,$ is uniformly equicontinuous
by the metric $\rho$ in $\overline{{D_0}_P}$ defined
by~(\ref{eq5M}). In other words, for any $\varepsilon>0$ there is
$\delta=\delta(\varepsilon)>0$ such that $\rho(f_m(x),
f_m(y))<\varepsilon$ whenever $|x-y|<\delta$ and for every $m\in
{\Bbb N}.$

Moreover, there is a subsequence $f_{m_k},$ $k=1,2,\dots ,$ which
converges to $f$ uniformly by the metric $\rho.$ In this case, $f$
has a continuous boundary extension $f:\overline{D}\rightarrow
\overline{D_0}_P$ and, besides that, for any $x_0\in
\partial D$ there is $P_0:=f(x_0)\in E_{D_0}=\overline{D_0}_P\setminus D_0$ such that, for any
$\varepsilon>0$ there is $\delta=\delta(\varepsilon)>0$ and
$M=M(\varepsilon)\in {\Bbb N}$ such that $\rho(f_{m_k}(x),
P_0)<\varepsilon$ for all $x\in B(x_0, \delta)\cap D$ and
$k\geqslant M_0.$ }
\end{theorem}

\medskip
\begin{corollary}\label{cor1}
{\it\, The statement of Theorem~\ref{th2} remains true if, instead
of the condition regarding the integrability of the function $Q$
over spheres with respect to some set $E_1$ is replaced by a simpler
condition: $Q\in L^1(D_0).$}
\end{corollary}

\medskip
In some situations, the study of the class of mappings given below
may be more important than the previous one. Moreover, in fact, the
conclusion of Theorem~\ref{th2} may be obtained from another (more
general) result given below.

\medskip
Given $\delta>0,$ domains $D, D_0\subset {\Bbb R}^n,$ $n\geqslant
2,$ a compactum $E\subset D_0,$ a Lebesgue measurable function
$Q:{\Bbb R}^n\rightarrow [0, \infty],$ $Q(y)\equiv 0$ outside $D_0,$
we denote by $\frak{R}^{E}_{Q, \delta}(D, D_0)$ is a class of all
open, discrete and closed mappings $f:D\rightarrow D_0$ of a domain
$D$ onto some domain $f(D),$ $E\subset f(D)\subset D_0,$ such that
(\ref{eq2*A})--(\ref{eqA2}) hold for every $y_0\in \overline{D_0}$
and, in addition, $f^{\,-1}(E)$ is a continuum with
$h(f^{\,-1}(E))\geqslant \delta.$

\medskip
\begin{theorem}\label{th3}
{\it\, Let $\delta>0,$ let $D, D_0\subset {\Bbb R}^n$ be domains,
$n\geqslant 2,$ let $E$ be a compact set in $D_0$ and let $Q:{\Bbb
R}^n\rightarrow [0, \infty],$ $Q(y)\equiv 0$ outside $D_0,$ be a
Lebesgue measurable function. Let $f_m\in \frak{R}^{E}_{Q,
\delta}(D, D_0),$ $m=1,2,\ldots .$ Assume that, all the conditions
1)--5) of Theorem~\ref{th2} are fulfilled for $f_m,$ $m=1,2,\ldots
.$ Then the conclusion of Theorem~\ref{th2} holds for $f_m\in
\frak{R}^{E}_{Q, \delta}(D, D_0),$ $m=1,2,\ldots .$ }
\end{theorem}

\medskip
\begin{corollary}\label{cor2}
{\it\, The statement of Theorem~\ref{th3} remains true if, instead
of the condition regarding the integrability of the function $Q$
over spheres with respect to some set $E_1$ is replaced by a simpler
condition: $Q\in L^1(D_0).$}
\end{corollary}

\section{Preliminaries}
As usual, we use the notation
\begin{equation}\label{eq1_A_4}
C(f, x):=\{y\in \overline{{\Bbb R}^n}:\exists\,x_k\in D:
x_k\rightarrow x, f(x_k) \rightarrow y, k\rightarrow\infty\}\,.
\end{equation}
A mapping $f$ between domains $D$ and $D^{\,\prime}$ is called {\it
closed} if $f(E)$ is closed in $D^{\,\prime}$ for any closed set
$E\subset D$ (see, e.g., \cite[Section~3]{Vu}). Any open discrete
closed mapping is boundary preserving, i.e.  $C(f, \partial
D)\subset
\partial D^{\,\prime},$ where
\begin{equation}\label{eq1_A_5}
C(f, \partial D)=\bigcup\limits_{x\in \partial D}C(f, x)
\end{equation}
(see e.g.~\cite[Theorem~3.3]{Vu}). The following statement holds.

\medskip
\begin{lemma}\label{lem3} {\bf(V\"{a}is\"{a}l\"{a}'s lemma on the weak flatness of inner points).}
{\sl\, Let $n\geqslant 2 $, let $D$ be a domain in $\overline{{\Bbb
R}^n},$ and let $x_0\in D.$ Then for each $P>0$ and each
neighborhood $U$ of point $x_0$ there is a neighborhood $V\subset U$
of the same point such that $M(\Gamma(E, F, D))> P$ for any continua
$E, F \subset D $ intersecting $\partial U$ and $\partial V.$}
\end{lemma}

\medskip
The proof of Lemma~\ref{lem3} is essentially given by
V\"{a}is\"{a}l\"{a} in~\cite[(10.11)]{Va}, however, we have also
given a formal proof, see \cite[Lemma~2.2]{SevSkv$_2$}.~$\Box$

\medskip
The following statement holds, see, e.g.,
\cite[Theorem~1.I.5.46]{Ku}).

\medskip
\begin{proposition}\label{pr2}
{\it\, Let $A$ be a set in a topological space $X.$ If the set $C$
is connected and $C\cap A\ne \varnothing\ne C\setminus A,$ then
$C\cap
\partial A\ne\varnothing.$}
\end{proposition}

\medskip
A path $\alpha: [a,\,b)\rightarrow D$ is called a {\it total
$f$-lifting} of $\beta$ starting at $x,$ if $(1)\quad
\alpha(a)=x\,;$ $(2)\quad (f\circ\alpha)(t)=\beta(t)$ for any $t\in
[a, b).$ The following statement holds, see e.g.
\cite[Lemma~3.7]{Vu}.

\medskip
\begin{proposition}\label{pr4_a}
Let $f:D \rightarrow {\Bbb R}^n$ be a discrete open and closed
(boundary preserving) mapping, $\beta: [a,\,b)\rightarrow f(D)$ be a
path, and $x\in\,f^{-1}\left(\beta(a)\right).$ Then $\beta$ has a
total $f$-lif\-ting starting at $x.$
\end{proposition}

We set
$$
q_{y_0}(r)=\frac{1}{\omega_{n-1}r^{n-1}}\int\limits_{S(y_0,
r)}Q(y)\,d\mathcal{H}^{n-1}(y)\,,$$
and $\omega_{n-1}$ denotes the area of the unit sphere ${\Bbb
S}^{n-1}$ in ${\Bbb R}^n.$ Following \cite{MRSY}, we say that a
function ${\varphi}:D\rightarrow{\Bbb R}$ has a {\it finite mean
oscillation} at a point $x_0\in D,$ write $\varphi\in FMO(x_0),$ if
$$\limsup\limits_{\varepsilon\rightarrow
0}\frac{1}{\Omega_n\varepsilon^n}\int\limits_{B( x_0,\,\varepsilon)}
|{\varphi}(x)-\overline{{\varphi}}_{\varepsilon}|\ dm(x)<\infty\,,
$$
where $\overline{{\varphi}}_{\varepsilon}=\frac{1}
{\Omega_n\varepsilon^n}\int\limits_{B(x_0,\,\varepsilon)}
{\varphi}(x) \,dm(x)$ and $\Omega_n$ is the volume of the unit ball
${\Bbb B}^n$ in ${\Bbb R}^n.$
We also say that a function ${\varphi}:D\rightarrow{\Bbb R}$ has a
finite mean oscillation at $A\subset \overline{D},$ write
${\varphi}\in FMO(A),$ if ${\varphi}$ has a finite mean oscillation
at any point $x_0\in A.$

\medskip
Given a continuum $E\subset D,$ $\delta>0$ and a Lebesgue measurable
function $Q:{\Bbb R}^n\rightarrow [0, \infty]$ we denote by
$\frak{F}_{E, \delta}(D)$ the family of all open discrete mappings
$f:D\rightarrow {\Bbb R}^n,$ $n\geqslant 2,$ satisfying
relations~(\ref{eq2*A})--(\ref{eqA2}) at any point $y_0\in
\overline{{\Bbb R}^n}$ such that $h(f(E))\geqslant \delta.$ The
following statement holds (see \cite[Theorem~1.1]{ST}).

\medskip
\begin{proposition}\label{pr3}
{\it Let $D$ be a domain in ${\Bbb R}^n,$ $n\geqslant 2,$ and let
$B(x_0, \varepsilon_1)\subset D$ for some $\varepsilon_1>0.$

\medskip
Assume that, $Q\in L^1({\Bbb R}^n)$ and, in addition, one of the
following conditions hold:

\medskip
1) $Q\in FMO(\overline{{\Bbb R}^n});$

\medskip
2) for any $y_0\in \overline{{\Bbb R}^n}$ there is $\delta(y_0)>0$
such that $\int\limits_{\varepsilon}^{\delta(y_0)}
\frac{dt}{tq_{y_0}^{\frac{1}{n-1}}(t)}<\infty$ for every
$\varepsilon\in (0, \delta(y_0))$ and
%
$$
\int\limits_{0}^{\delta(y_0)}
\frac{dt}{tq_{y_0}^{\frac{1}{n-1}}(t)}=\infty\,.
$$
Then there is $r_0>0,$ which does not depend on $f,$ such that
$$f(B(x_0, \varepsilon_1))\supset B_h(f(x_0), r_0)\qquad \forall\,\,f\in \frak{F}_{E,
\delta}(D)\,,$$
where $B_h(f(x_0), r_0)=\{w\in \overline{{\Bbb R}^n}: h(w,
f(x_0))<r_0\}.$ }
\end{proposition}

\section{Main Lemma}

The following lemma holds, cf. \cite[Lemma~4.1]{SevSkv$_2$},
\cite[Lemma~2.13]{ISS}.

\medskip
\begin{lemma}\label{lem1}
{\it\, Let $\delta>0,$ let $D, D_0\subset {\Bbb R}^n$ be domains,
$n\geqslant 2,$ let $E$ be a compact set in $D_0$ and let $Q:{\Bbb
R}^n\rightarrow [0, \infty],$ $Q(y)\equiv 0$ outside $D_0,$ be a
Lebesgue measurable function. Assume that, no connected component of
the boundary of the domain $D$ degenerates into a point and, besides
that, $f_{m}(D)$ converge to $D_0$ as its kernel for $f_m\in
\frak{F}^{E}_{Q, \delta}(D, D_0),$ $m=1,2,\ldots .$ Assume that, for
each point $y_0\in \overline{D_0}$ and for every
$0<r_1<r_2<r_0:=\sup\limits_{y\in D_0}|y-y_0|$ there is a set
$E_1\subset[r_1, r_2]$ of a positive linear Lebesgue measure such
that the function $Q$ is integrable with respect to
$\mathcal{H}^{n-1}$ over the spheres $S(y_0, r)$ for every $r\in
E_1.$

\medskip
If $E_*$ is some another compactum in $D_0$ with $E_*\subset
f_m(D)\subset D_0$ for every $m\in {\Bbb N},$ then there exists
$\delta_*>0$ such that $h(f_m^{\,-1}(E_*),
\partial D)\geqslant \delta_*>0$ for every $m\in {\Bbb N}.$
}
\end{lemma}

\medskip
\begin{proof} Since $D$ is a domain in ${\Bbb R}^n,$ $\partial
D\ne\varnothing.$ Thus, the quantity $h(f_m^{\,-1}(E_*),
\partial D)$ is well-defined.

\medskip
Let us prove Lemma~\ref{lem1} by the contradiction. Assume that, the
conclusion of the lemma is not true. Then for each $k\in {\Bbb N}$
there is some number $m_k\in {\Bbb N}$ such that
$h(f^{\,-1}_{m_k}(E_*),
\partial D)<1/k.$ Since $f_m$ is open, discrete and closed, the set
$f^{\,-1}_{m}(E_*)$ is a compactum in $D$ for each $m=1,2,\ldots$
(see e.g. \cite[Theorem~3.3]{Vu}). Therefore, we may consider that
the sequence $m_k$ is increasing by $k=1,2,\ldots .$ Renumbering the
elements $f_{m_k},$ if required, we may consider that the latter
holds for the same sequence $f_m,$ i.e., for each $m\in {\Bbb N}$
there is some number $m\in {\Bbb N}$ such that $h(f^{\,-1}_{m}(E_*),
\partial D)<1/m.$

\medskip
Since $f^{\,-1}_{m}(E_*)$ is a compactum in $D$ for each
$m=1,2,\ldots,$ we have that $h(f^{\,-1}_{m}(E_*),
\partial D)=h(x_m, y_m)<\frac{1}{m}$ for some
$x_m\in f^{\,-1}_{m}(E_*)$ and $y_m\in \partial D.$
Since $\partial D$ is a compact set, we may assume that
$y_m\rightarrow y_0$ as $m\rightarrow\infty$ for some $y_0\in
\partial D.$ Now, $x_m\rightarrow y_0$ as $m\rightarrow\infty,$ as
well.

\medskip
Let $w_m=f_m(x_m)\in E_*.$ We may consider that $w_m$ converges to
some point $w_0$ in $E_*$ as $m\rightarrow\infty.$ Observe that,
$w_m\in B(w_0, r_0)\subset D$ for some $r_0>0$ and sufficiently
large $m=1,2\ldots.$ Let $z_0$ be some point in $E.$ Let us join the
points $w_0$ and $z_0$ by some path $\gamma:[1, 2]\rightarrow D_0,$
$\gamma(1)=w_0,$ $\gamma(2)=z_0,$ and let $\gamma_m:[0,
1]\rightarrow D_0$ be a segment $\gamma_m(t)=w_m+(w_0-w_m)t,$ $t\in
[0, 1].$ Since $D_0$ is a kernel of $f_m(D)$ and
$K_0:=|\gamma|\cup\overline{B(w_0, r_0)}$ is a compactum in $D_0,$
we may assume that $K_0\subset f_m(D)$ for all $m=1,2,\ldots .$ Set
$$E_m(t)=\begin{cases}\gamma_m(t)\,,& t\in [0, 1]\,\\
\gamma(t)\,,& t\in [1, 2]\end{cases}\,.$$
Now, $|E_m|$ is a compactum in $K_0$ for all $m=1,2,\ldots.$ Let
$A_m:[0, 2)\rightarrow D$ be a whole lifting of $E_m,$ $|E_m|\subset
f_m(D),$ starting at $x_m$ (it exists by Proposition~\ref{pr4_a}).

\medskip
Observe that, no path $A_m(t),$ $A_m:[0, 2)\rightarrow D,$ cannot
tend to the boundary of the domain $D$ as $t\rightarrow 2-0,$
because $C(f_m, \partial D)\subset \partial f_m(D)$ whenever $f_m$
is an open, discrete and closed mapping (see e.g.
\cite[Theorem~3.3]{Vu}). Let us to prove that $A_m$ has a limit as
$t\rightarrow 2-0.$ Assume the contrary, i.e., there is $m\in {\Bbb
N}$ such that a path $A_m$ has no a limit as $t\rightarrow 2-0.$
Then the cluster set $C(A_m(t), 2)$ is a continuum in $D.$ Since the
mapping $f_m$ is continuous in $D$ for any $m\in {\Bbb N},$ we
obtain that $f_m\equiv const$ on $C(A_m(t), 2),$ which contradicts
the discreteness of $f_m.$

\medskip
Now, let $A_m:[0, 2]\rightarrow D$ be a continuous extension of
$A_m$ at $t=2$ (here we preserve the notion of $A_m$ for
simplicity). Observe that, $A_m(2):=z_m\in f_m^{\,-1}(E)$ by the
definition of $f_m^{\,-1}(E).$

\medskip
Let $E_0$ be a component of $\partial D$ containing $y_0.$ Since by
the assumptions of the lemma, all components of $\partial D$ are
non-degenerate, there exists $r>0$ such that $h(E_0)\geqslant r.$
Put $P>0$ and $U=B_h(y_0, R_0)=\{y\in \overline{{\Bbb R}^n}: h(y,
y_0)<R_0\},$ where $2R_0:=\min\{r/2, \delta/2\}$ and $\delta$ is a
number in the definition of the class $\frak{F}^{E}_{Q, \delta}(D,
D_0)$. Observe that $A_m\cap U\ne\varnothing\ne A_m\setminus U$ for
sufficiently large $m\in{\Bbb N},$ since $x_m\rightarrow y_0$ as
$m\rightarrow \infty,$ $x_m\in A_m;$ besides that $h(A_m)\geqslant
\delta\geqslant 2R_0$ and $h(U)\leqslant 2R_0.$ Since $A_m$ is a
continuum, $A_m\cap
\partial U\ne\varnothing$ by Proposition~\ref{pr2}. Similarly,
$E_0\cap U\ne\varnothing\ne E_0\setminus U$ for sufficiently large
$m\in{\Bbb N},$ since $h(E_0)\geqslant r> 2R_0$ and $h(U)\leqslant
2R_0.$ Since $E_0$ is a continuum, $E_0\cap
\partial U\ne\varnothing$ by Proposition~\ref{pr2}.
By the proving above,
\begin{equation}\label{eq8}
A_m\cap
\partial U\ne\varnothing\ne E_0\cap
\partial U\,.
\end{equation}
By Lemma~\ref{lem3} there is $V\subset U,$ $V$ is a neighborhood of
$y_0,$ such that
\begin{equation}\label{eq9}
M(\Gamma(E, F, \overline{{\Bbb R}^n}))>P
\end{equation}
for any continua $E, F\subset \overline{{\Bbb R}^n}$ with $E\cap
\partial U\ne\varnothing\ne E\cap \partial V$ and $F\cap \partial
U\ne\varnothing\ne F\cap \partial V.$
Arguing similarly to above, we may prove that
$$
A_m\cap
\partial V\ne\varnothing\ne E_0\cap
\partial V
$$
for sufficiently large $m\in {\Bbb N}.$ Thus, by~(\ref{eq9})
\begin{equation}\label{eq9B}
M(\Gamma(A_m, E_0,\overline{{\Bbb R}^n}))>P
\end{equation}
for sufficiently large $m=1,2,\ldots ,$ see Figure~\ref{fig1}.
\begin{figure}[h]
\centerline{\includegraphics[scale=0.5]{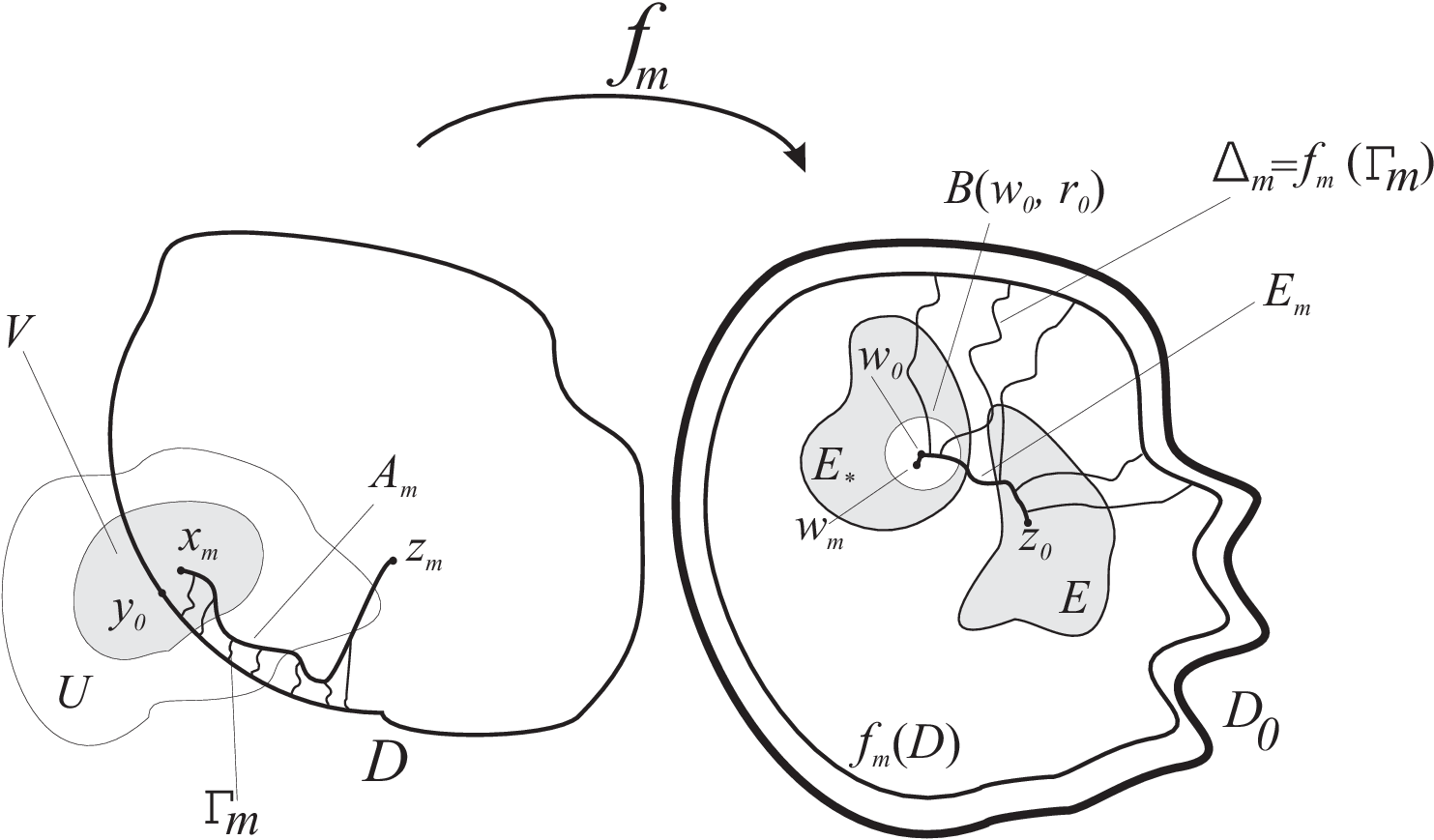}} \caption{To proof
of Lemma~\ref{lem1}}\label{fig1}
\end{figure}
Let $\gamma:[0, 1]\rightarrow \overline{{\Bbb R}^n}$ be a path in
$\Gamma(A_m, E_0, \overline{{\Bbb R}^n}),$ i.e., $\gamma(0)\in A_m,$
$\gamma(1)\in E_0$ and $\gamma(t)\in \overline{{\Bbb R}^n}$ for
$t\in (0, 1).$ Let $t_m=\sup\limits_{\gamma(t)\in D}t$ and let
$\alpha_m(t)=\gamma|_{[0, t_m)}(t).$ Let $\Gamma_m$ consists of all
such paths $\alpha_m,$ now $\Gamma(A_m, E_0, \overline{{\Bbb
R}^n})>\Gamma_m$ and by the minorization principle of the modulus
(see~\cite[Theorem~1]{Fu})
\begin{equation}\label{eq12}
M(\Gamma_m)\geqslant M(\Gamma(A_m, E_0, \overline{{\Bbb R}^n}))\,.
\end{equation}
Combining~(\ref{eq9B}) and~(\ref{eq12}), we obtain that
\begin{equation}\label{eq9C}
M(\Gamma_m)>P
\end{equation}
for sufficiently large $m=1,2,\ldots .$

\medskip
We now prove that the relation~(\ref{eq9C}) contradicts the
definition of $f_m$ in~(\ref{eq2*A})--(\ref{eqA2}). Since $f_m$ is a
closed mapping, it preserves the boundary (see
\cite[Theorem~3.3]{Vu}), so that $C(f_m,
\partial D)\subset \partial f_m(D).$ Thus, $C(\beta_m(t), t_m)\subset \partial f_m(D),$
$\beta_m(t):=f_m(\alpha_m(t))=f_m(\gamma|_{[0, t_m)}(t)),$
$\gamma\in \Gamma(A_m, E_0, \overline{{\Bbb R}^n})$ and
$t_m=\sup\limits_{\gamma(t)\in D}t.$

\medskip
Observe that, ${\rm dist\,}(K_0, \partial f_m(D))\geqslant
\varepsilon>0$ for some $\varepsilon>0$ and all sufficiently large
$m=1,2,\ldots.$ Otherwise, there is a sequence of numbers
$\varepsilon_m>0,$ $m\in {\Bbb N},$ and points $p_m\in K_0,$ $q_m\in
\partial f_m(D)$ such that $|p_m-q_m|<1/m.$ Since $K_0$ is a
compactum in $D_0,$ we may consider that $p_m\rightarrow p_0$ as
$m\rightarrow\infty$ for some $p_0\in K_0.$ Since $K_0\subset D_0$
by the construction, besides that, the sequence $f_m(D)$ converges
to $D_0$ as the kernel, $p_0\subset f_m(D)$ for sufficiently large
$m\in {\Bbb N.}$ Moreover, let $\widetilde{\varepsilon}:={\rm
dist}\,(p_0,\partial D_0),$ now $p_0, p_m\in B(p_0,
\widetilde{\varepsilon}/2)\subset D_0$ for sufficiently large $m\in
{\Bbb N}$ and $B(p_0, \widetilde{\varepsilon}/2)$ contains in
$f_m(D)$ for sufficiently large $m\in {\Bbb N}.$ But, at the same
time, $q_m\rightarrow p_0$ as $m\rightarrow\infty$ and,
consequently, $\partial f_m(D)\cap \overline{B(p_0,
\widetilde{\varepsilon}/2)}=\varnothing$ for sufficiently large
$m\in {\Bbb N}.$ The latter contradiction disproves the assumption
made above.

\medskip
Now, ${\rm dist\,}(K_0, \partial f_m(D))\geqslant \varepsilon>0$ for
some $\varepsilon>0$ and all sufficiently large $m=1,2,\ldots.$ We
cover the continuum $K_0$ with balls $B(x, \varepsilon/4),$ $x\in
A.$ Since $K_0$ is a compact set, we may assume that $K_0\subset
\bigcup\limits_{i=1}^{M_0}B(x_i, \varepsilon/4),$ $x_i\in K_0,$
$i=1,2,\ldots, M_0,$ $1\leqslant M_0<\infty.$ By the definition,
$M_0$ depends only on $K_0,$ in particular, $M_0$ does not depend on
$m.$
Note that
\begin{equation}\label{eq6D}
\Gamma_m=\bigcup\limits_{i=1}^{M_0}\Gamma_{mi}\,,
\end{equation}
where $\Gamma_{mi}$ consists of all paths $\gamma:[0,
t_m)\rightarrow D$ in $\Gamma_m$ such that $f_m(\gamma(0))\in B(y_i,
\varepsilon/4).$ We now show that
\begin{equation}\label{eq7C}
f_m(\Gamma_{mi})>\Gamma(S(y_i, \varepsilon/4), S(y_i,
\varepsilon/2), A(y_i, \varepsilon/4, \varepsilon/2))\,.
\end{equation}
Indeed, let $\widetilde{\gamma}\in f_m(\Gamma_{mi}),$ now there is
$\gamma\in \Gamma_{mi}$ such that $\widetilde{\gamma}=f_m(\gamma).$
If $\gamma\in \Gamma_{mi},$ then $\gamma:[0, t_m)\rightarrow D,$
$\gamma\in \Gamma_m,$ $f_m(\gamma(0))\in B(y_i, \varepsilon/4)$ and
$\gamma(t)\rightarrow
\partial D$ as $t\rightarrow t_m-0.$ By the comments mentioned
above, $\widetilde{\gamma}(t)=f_m(\gamma(t))\rightarrow
\partial f_m(D)$ as $t\rightarrow t_m-0.$ Now, by the definition of $\varepsilon,$
$|\widetilde{\gamma}|\cap B(y_i, \varepsilon/4)\ne\varnothing\ne
|\widetilde{\gamma}|\cap (D_0\setminus B(y_i, \varepsilon/4)).$
Therefore, by Proposition~\ref{pr2} there is $0<t_1<1$ such that
$\widetilde{\gamma}(t_1)\in S(y_i, \varepsilon/4).$ We may assume
that $\widetilde{\gamma}(t)\not\in B(y_i, \varepsilon/4)$ for
$t>t_1.$ Put $\gamma_1:=\widetilde{\gamma}|_{[t_1, 1]}.$ Similarly,
$|\gamma_1|\cap B(y_i, \varepsilon/2)\ne\varnothing\ne
|\gamma_1|\cap (D_0\setminus B(y_i, \varepsilon/2)).$ By
Proposition~\ref{pr2} there is $t_1<t_2<1$ with $\gamma_1(t_2)\in
S(y_i, \varepsilon/2).$ We may assume that $\gamma_1(t)\in B(y_i,
\varepsilon/2)$ for $t<t_2.$ Put $\gamma_2:=\gamma|_{[t_1, t_2]}.$
Then, the path $\gamma_2$ is a subpath of $\widetilde{\gamma},$
which belongs to the family $\Gamma(S(y_i, \varepsilon/4), S(y_i,
\varepsilon/2), A(y_i, \varepsilon/4, \varepsilon/2)).$ Thus, the
relation~(\ref{eq7C}) is established. By~(\ref{eq6D}) and
(\ref{eq7C}), we obtain that
\begin{equation}\label{eq5}
\Gamma_m>\bigcup\limits_{i=1}^{M_0}\Gamma_{f_m}(y_i, \varepsilon/4,
\varepsilon/2)\,.
\end{equation}
We set $\widetilde{Q}(y)=\max\{Q(y), 1\}$ and
$$\widetilde{q}_{y_i}(r)=\int\limits_{S(y_i,
r)}\widetilde{Q}(y)\,d\mathcal{A}\,.$$ Now,
$\widetilde{q}_{y_i}(r)\ne \infty$ for $r\in
E_1\subset[\varepsilon/4, \varepsilon/2],$ where $E_1$ is some set
of positive linear measure which exists by the assumptions of the
lemma. Set
$$I_i=I_i(y_i, \varepsilon/4, \varepsilon/2)=\int\limits_{\varepsilon/4}^{\varepsilon/2}\
\frac{dr}{r\widetilde{q}_{y_i}^{\frac{1}{n-1}}(r)}\,.$$
Observe that $I\ne 0,$ because $\widetilde{q}_{y_i}(r)\ne \infty$
for $r\in E_1\subset[\varepsilon/4, \varepsilon/2],$ where $E_1$ is
some set of positive linear measure. Besides that, by the direct
calculations we obtain that
$$I_i\leqslant \log\frac{r_2}{r_1}<\infty\,,\quad i=1,2, \ldots, M_0\,.$$
Now, we put
$$\eta_i(r)=\begin{cases}
\frac{1}{I_ir\widetilde{q}_{y_i}^{\frac{1}{n-1}}(r)}\,,&
r\in [\varepsilon/4, \varepsilon/2]\,,\\
0,& r\not\in [\varepsilon/4, \varepsilon/2]\,.
\end{cases}$$
Observe that, a function~$\eta_i$ satisfies the
condition~$\int\limits_{\varepsilon/4}^{\varepsilon/2}\eta_i(r)\,dr=1,$
therefore it may be substituted into the right side of the
inequality~(\ref{eq2*A}) with the corresponding values $f,$ $r_1$
and $r_2.$ Now, we obtain that
\begin{equation}\label{eq7B}
M(\Gamma_{f_m}(y_i, \varepsilon/4, \varepsilon/2))\leqslant
\int\limits_{A(y_i, \varepsilon/4, \varepsilon/2)}
\widetilde{Q}(y)\,\eta^n_i(|y- y_i|)\,dm(y)\,.\end{equation}
By the Fubini theorem we have that
$$\int\limits_{A(y_i, \varepsilon/4, \varepsilon/2)}
\widetilde{Q}(y)\,\eta^n_i(|y- y_i|)\,dm(y)=$$
\begin{equation}\label{eq7E}
=\int\limits_{\varepsilon/4}^{\varepsilon/2}\int\limits_{S(y_i,
r)}Q(y)\eta^n_i(|y- y_i|)\,d\mathcal{A}\,dr\,=
\end{equation}
$$=\frac{\omega_{n-1}}{I_i^n}\int\limits_{\varepsilon/4}^{\varepsilon/2}r^{n-1}
\widetilde{q}_{y_i}(r)\cdot
\frac{dr}{r^n\widetilde{q}^{\frac{n}{n-1}}_{y_i}(r)}=\frac{\omega_{n-1}}{I_i^{n-1}}\,,$$
where $\omega_{n-1}$ is the area of the unit sphere in ${\Bbb R}^n.$
Now, by~(\ref{eq7B}) and~(\ref{eq7E}) we obtain that
$$M(\Gamma_{f_m}(y_i, \varepsilon/4,
\varepsilon/2))\leqslant \frac{\omega_{n-1}}{I_i^{n-1}}\,,$$
whence from~(\ref{eq5}) we obtain that
\begin{equation}\label{eq7D}
M(\Gamma_m)\leqslant \sum\limits_{i=1}^{M_0}M(\Gamma_{f_m}(y_i,
\varepsilon/4, \varepsilon/2))\leqslant
\sum\limits_{i=1}^{M_0}\frac{\omega_{n-1}}{I_i^{n-1}}:=C_0\,, \quad
m=1,2,\ldots\,.
\end{equation}
Since $P$ in~(\ref{eq9B}) may be done arbitrary big, the
relations~(\ref{eq9B}) and~(\ref{eq7D}) contradict each other. This
completes the proof.~$\Box$
\end{proof}

\section{Proof of Theorem~\ref{th2}}

For the proof, we use some approaches developed in~\cite{Sev$_2$}.
Let us firstly prove that the family $f_m,$ $m=1,2,\ldots ,$ is
uniformly equicontinuous by the metric $\rho$ in
$\overline{{D_0}_P}$ which is defined in~(\ref{eq5M}). In other
words, we need to prove that, for any $\varepsilon>0$ there is
$\delta=\delta(\varepsilon)>0$ such that $\rho(f_m(x),
f_m(y))<\varepsilon$ whenever $|x-y|<\delta$ and for every $m\in
{\Bbb N}.$

\medskip
We will prove it by contradiction, namely, assume that there is
$\varepsilon>0$ such that, for any $k\in {\Bbb N}$ there is
$m=m_k\in {\Bbb N}$ and $x_k, y_k\in D$ such that
$|x_k-y_k|<\frac{1}{k},$ however, $\rho(f_{m_k}(x_k),
f_{m_k}(y_{k}))\geqslant\varepsilon.$ Since any $f_m$ has a
continuous extension to $\overline{D},$ we may assume that the
sequence $m_k$ is increasing by $k.$ Now, $f_{m_k}(D)$ converges to
$D_0$ as to its kernel as $k\rightarrow\infty,$ as well. By
resorting to renumbering, if necessary, we may assume that the
sequence $f_m$ itself satisfies the above condition, i.e.,
$|x_m-y_m|<\frac{1}{m},$ $m=1,2,\ldots,$ however,
\begin{equation}\label{eq1A}
\rho(f_m(x_m), f_m(y_m))>\varepsilon\,,\qquad m=1,2,\ldots\,.
\end{equation}
Due to the compactness of $\overline{{\Bbb R}^n},$ we may assume
that $x_m, y_m\rightarrow x_0$ as $m\rightarrow\infty.$ Besides
that, since the space $\overline{D_0}_P$ is compact, we may assume
that $f_m(x_m)$ and $f_m(y_m)$ converge so some $P_1\ne P_2$ as
$m\rightarrow\infty,$ $P_1, P_2\in \overline{D_0}_P.$

\medskip
Let $\sigma_m$ and $\sigma^{\,\prime}_m,$ $m=0,1,2,\ldots, $ be
sequences of cuts corresponding $P_1$ and $P_2,$ respectively. Let
$\sigma_m,$ $m=0,1,2,\ldots, $ lie on the spheres $S(z_0, r_m)$
centered at some point $z_0\in
\partial D_0,$ where $r_m\rightarrow 0$ as $m\rightarrow\infty$
(such a sequence $\sigma_m$ exists by~\cite[Lemma~3.1]{IS},
cf.~\cite{KR}). If at least one of the points $P_1$ (or $P_2$) are
inner points of $D_0,$ then $\sigma_m$ are radii of balls centered
at $P_1$ (or~$P_2$).

Let $d_m$ and $g_m,$ $m=0,1,2,\ldots, $ be sequences of domains in
$D_0,$ which correspond to cuts $\sigma_m$ and
$\sigma^{\,\prime}_m,$ respectively. If at least one of the points
$P_1$ (or $P_2$) are inner points of $D_0,$ then $d_m=B(P_1,
\sigma_m)$ (or $g_m=B(P_2, \kappa_m)$), where $\sigma_m$ are radii
of balls centered at $P_1$ (or $P_2$). Since $(\overline{D_0}_P,
\rho)$ is a metric space, we may consider that $d_m$ and $g_m$ are
disjoint for every $m=0,1,2,\ldots ,$ in particular,
\begin{equation}\label{eq4}
d_0\cap g_0=\varnothing\,.
\end{equation}
Since $f_m(D)$ is a regular sequence of domains with a respect to
$D_0,$ there exists a sequence of cuts $\varsigma_{m}$ and domains
$d^{\,\prime}_m,$ equivalent to $\sigma_m$ and $d_m,$ respectively,
such that the following condition holds: given $d^{\,\prime}_k,$
$k=1,2,\ldots,$ there is $M_1=M_1(k)$ such that $d^{\,\prime}_k\cap
f_m(D)$ is a non-empty connected set for every $m\geqslant M_1(k).$
Similarly, we may consider that there exists a sequence of cuts
$\varsigma^{\,\prime}_{m}$ and domains $g^{\,\prime}_m,$ equivalent
to $\sigma^{\,\prime}_m$ and $g_m,$ respectively, such that the
following condition holds: given $g^{\,\prime}_k,$ $k=1,2,\ldots,$
there is $M_2=M_2(k)$ such that $g^{\,\prime}_k\cap f_m(D)$ is a
non-empty connected set for every $m\geqslant M_2(k).$ Now, we will
consider that $d^{\prime}_k\subset d_k$ and $g^{\prime}_k\subset
g_k$ for any $k=1,2,\ldots .$ Since $f_m(x_m)$ converge to $P_1$ as
$m\rightarrow\infty,$ we may consider that $f_m(x_m)\in
d^{\prime}_m,$ $m=1,2,\ldots.$ Similarly, we may consider that
$f_m(y_m)\in g^{\,\prime}_m,$ $m=1,2,\ldots.$

\begin{figure}[h]
\centerline{\includegraphics[scale=0.55]{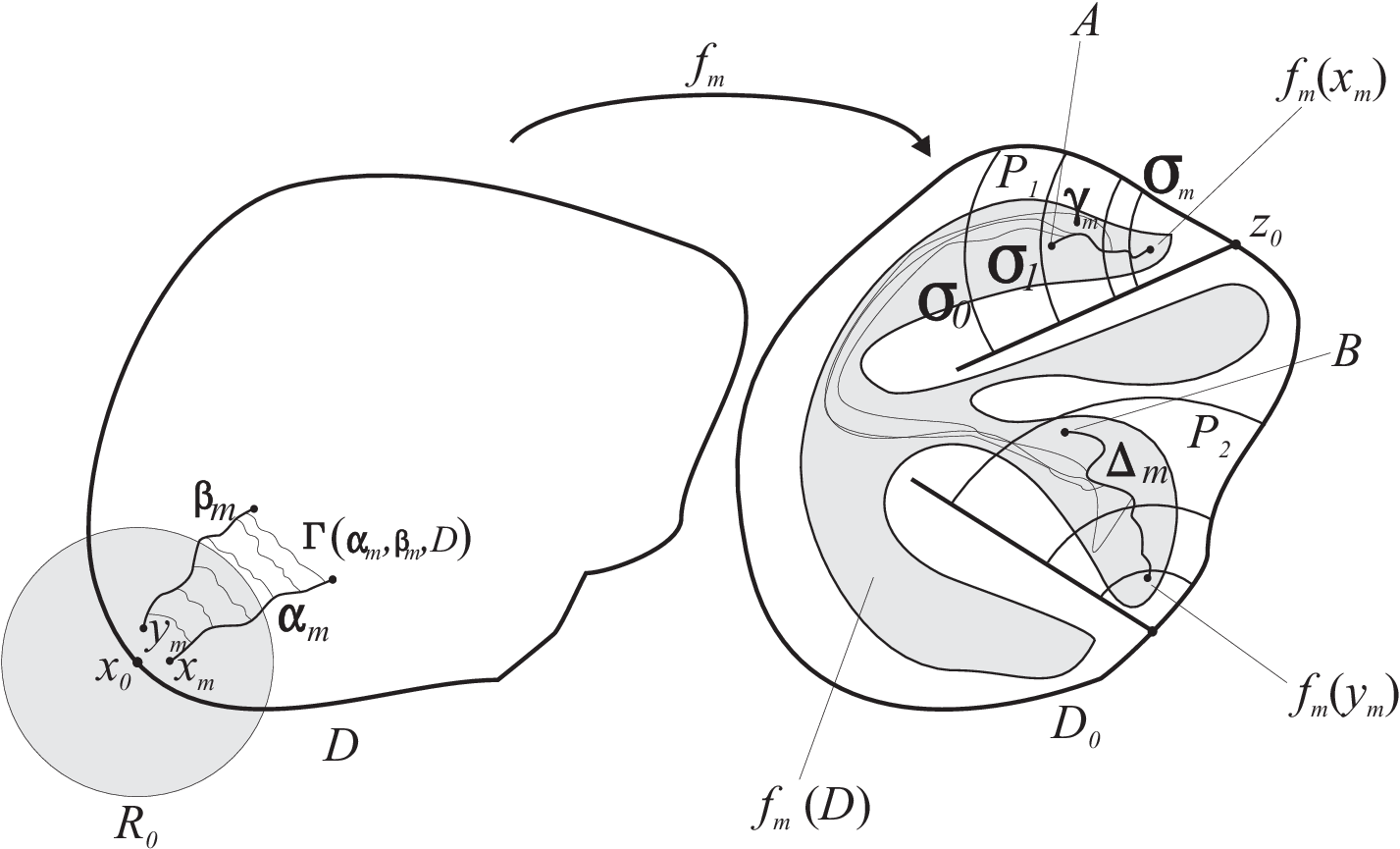}}
\caption{To the proof of Theorem~\ref{th2}}\label{fig1A}
\end{figure}

\medskip
Put $A\in d^{\,\prime}_1$ and $B\in g^{\,\prime}_1.$ Since
$d^{\,\prime}_1, g^{\,\prime}_1\subset D_0$ and $f_m(D)$ converges
to $D_0$ as its kernel, then $A, B\in f_m(D)$ for $m\geqslant M_3$
and some sufficiently large $M_3\in {\Bbb N}.$ In addition, due to
the mentioned above, $f_m(x_m)\in d^{\prime}_1,$ $f_m(y_m)\in
g^{\,\prime}_1$ for $m=1,2,\ldots.$ Besides that, by the saying
above, $d^{\,\prime}_1\cap f_m(D)$ and $g^{\,\prime}_1\cap f_m(D)$
are connected sets for $m\geqslant \max\{M_1(1), M_2(1)\}.$ Now, let
$m\geqslant \max\{M_1(1), M_2(1), M_3\}.$ For such $m,$ we have that
$A, f_m(x_m)\in d^{\,\prime}_1\cap f_m(D)$ and $B, f_m(y_m)\in
g^{\,\prime}_1\cap f_m(D).$ Let $\gamma_m:[0, 1]\rightarrow D_0$ be
a path joining $f_m(x_m)$ and $A$ in $d^{\,\prime}_1\cap f_m(D),$
i.e., $\gamma_m(0)=f_m(x_m),$ $\gamma_m(1)=A$ and $\gamma_m(t)\in
d^{\,\prime}_1\cap f_m(D)$ for all $t\in (0, 1).$ Similarly, let
$\Delta_m:[0, 1]\rightarrow D_0$ be a path joining $f_m(y_m)$ and
$B$ in $g^{\,\prime}_1\cap f_m(D),$ i.e., $\Delta_m(0)=f_m(y_m),$
$\Delta_m(1)=B$ and $\Delta_m(t)\in g^{\,\prime}_1\cap f_m(D)$ for
all $t\in (0, 1).$ Let $\alpha_m$ and $\beta_m$ be total
$f_m$-liftings of $\gamma_m$ and $\Delta_m$ in $D$ starting at
points $x_m$ and $y_m,$ respectively (such liftings exist by
Proposition~\ref{pr4_a}), see Figure~\ref{fig1A}).
Now, he paths $\alpha_m$ and $\beta_m$ are starting at the points
$x_m$ and $y_m,$ respectively, and are ending at some points $C_m\in
f^{\,-1}_m(A)$ and $D_m\in f^{\,-1}_m(A).$ By Lemma~\ref{lem1}
sequently applied for $E_*:=A$ and $E_*=B,$ we have that
$h(f_m^{\,-1}(A),
\partial D)\geqslant \delta_*>0$ and $h(f_m^{\,-1}(B),
\partial D)\geqslant \delta_*>0$ for every $m\in {\Bbb N}.$ Then there is
$R_0>0$ such that $C_m, D_m\in D\setminus B(x_0, R_0)$ for all
$m=1,2,\ldots .$ Since the boundary of $D$ is weakly flat, for any
$P>0$ there is $m=m_P\geqslant 1$ such that
\begin{equation}\label{eq7}
M(\Gamma(|\alpha_m|, |\beta_m|, D))>P\qquad\forall\,\,m\geqslant
m_P\,.
\end{equation}
Observe that~(\ref{eq7}) holds also in the case, when $x_0$ is the
inner point of $D,$ because the inner points of any domain are
``weakly flat'', see Lemma~\ref{lem3}.

\medskip
We now show that~(\ref{eq7}) contradicts with~(\ref{eq2*A}). Indeed,
let $\gamma\in \Gamma(|\alpha_m|, |\beta_m|, D).$ Then $\gamma:[0,
1]\rightarrow D,$ $\gamma(0)\in |\alpha_m|$ and $\gamma(1)\in
|\beta_m|.$ In particular, $f_m(\gamma(0))\in |\gamma_m|$ and
$f_m(\gamma(1))\in |\Delta_m|.$ Now, by~(\ref{eq4}) and~(\ref{eq7})
it follows that $|f_m(\gamma)|\cap d_1\ne\varnothing \ne
|f_m(\gamma)|\cap(D_0\setminus d_1)$ for sufficiently large $m\in
{\Bbb N}.$ By~\cite[Theorem~1.I.5.46]{Ku} $|f_m(\gamma)|\cap
\partial d_1\ne\varnothing,$ i.e., $|f_m(\gamma)|\cap S(z_0,
r_1)\ne\varnothing,$ because $\partial d_1\cap D_0\subset
\sigma_1\subset S(z_0, r_1)$ by the definition of the cut
$\sigma_1.$ Let $t_1\in (0,1)$ be such that $f_m(\gamma(t_1))\in
S(z_0, r_1)$ and $f_m(\gamma)|_1:=f_m(\gamma)|_{[t_1, 1]}.$ Without
loss of generality, we may consider $f_m(\gamma)|_1\subset {\Bbb
R}^n\setminus d_1.$ Arguing similarly for $f_m(\gamma)|_1,$ we may
consider a point $t_2\in (t_1, 1)$ such that $f_m(\gamma(t_2))\in
S(z_0, r_0).$ We set $f(\gamma)|_2:=f(\gamma)|_{[t_1, t_2]}.$ Then
$f_m(\gamma)|_2$ is a subpath of $f_m(\gamma)$ and, in addition,
$f_m(\gamma)|_2\in \Gamma(S(z_0, r_1), S(z_0, r_0), D_0).$ Without
loss of generality we may assume that $f_m(\gamma)|_2\subset A(z_0,
r_1, r_0).$ Thus, $\Gamma(|\alpha_m|, |\beta_m|,
D)>\Gamma_{f_m}(z_0, r_1, r_0).$
From the latter relation, by the minorization principle of the
modulus (see \cite[Theorem~1(c)]{Fu})
\begin{equation}\label{eq5_1}
M(\Gamma(|\alpha_m|, |\beta_m|, D))\leqslant M(\Gamma_{f_m}(z_0,
r_1, r_0))\,.
\end{equation}
Now we argue as under the latter part of the proof of
Lemma~\ref{lem1}. We set $\widetilde{Q}(y)=\max\{Q(y), 1\}$ and
$$\widetilde{q}_{z_0}(r)=\int\limits_{S(z_0,
r)}\widetilde{Q}(y)\,d\mathcal{A}\,.$$
Now, $\widetilde{q}_{z_0}(r)\ne \infty$ for $r\in E_1\subset[r_1,
r_0],$ where $E_1$ is some set of positive linear measure which
exists by the assumption~1) of the theorem. Set
$$I=I(z_0, r_1, r_0)=\int\limits_{r_1}^{r_0}\
\frac{dr}{r\widetilde{q}_{y_i}^{\frac{1}{n-1}}(r)}\,.$$
Due to the mentioned above, $I\ne 0.$ Besides that, by the direct
calculations we obtain that
$$I\leqslant \log\frac{r_1}{r_0}<\infty\,.$$
Now, we put
$$\eta(r)=\begin{cases}
\frac{1}{Ir\widetilde{q}_{z_0}^{\frac{1}{n-1}}(r)}\,,& r\in [r_1,
r_0]\,,\\
0,& r\not\in [r_1, r_0]\,.
\end{cases}$$
Observe that, a function~$\eta$ satisfies the
condition~$\int\limits_{r_1}^{r_0}\eta(r)\,dr=1,$ therefore it may
be substituted into the right side of the inequality~(\ref{eq2*A})
with the corresponding values $f,$ $r_1$ and $r_2.$ Now, we obtain
that
\begin{equation}\label{eq7F}
M(\Gamma_{f_m}(z_0, r_1, r_0))\leqslant \int\limits_{A(z_0, r_1,
r_0)} \widetilde{Q}(y)\,\eta^n(|y- z_0|)\,dm(y)\,.\end{equation}
By the Fubini theorem we have that
$$\int\limits_{A(z_0, r_1,
r_0)} \widetilde{Q}(y)\,\eta^n_i(|y- z_0|)\,dm(y)=$$
\begin{equation}\label{eq7G}
=\int\limits_{r_1}^{r_0}\int\limits_{S(z_0, r)}Q(y)\eta^n_i(|y-
z_0|)\,d\mathcal{A}\,dr\,=
\end{equation}
$$=\frac{\omega_{n-1}}{I^n}\int\limits_{r_1}^{r_0}r^{n-1}
\widetilde{q}_{z_0}(r)\cdot
\frac{dr}{r^n\widetilde{q}^{\frac{n}{n-1}}_{z_0}(r)}=\frac{\omega_{n-1}}{I^{n-1}}\,,$$
where $\omega_{n-1}$ is the area of the unit sphere in ${\Bbb R}^n.$
Now, by~(\ref{eq7F}) and~(\ref{eq7G}) we obtain that
$$M(\Gamma_{f_m}(z_0, r_1, r_0))\leqslant \frac{\omega_{n-1}}{I^{n-1}}\,,$$
whence from~(\ref{eq5_1}) we obtain that
\begin{equation}\label{eq7H}
M(\Gamma(|\alpha_m|, |\beta_m|, D))\leqslant
\frac{\omega_{n-1}}{I^{n-1}}:=C_0\,, \quad m=1,2,\ldots\,.
\end{equation}
Since $P$ in~(\ref{eq7}) may be done arbitrary big, the
relations~(\ref{eq7}) and~(\ref{eq7H}) contradict each other. The
obtained contradiction proves the uniform equicontinuity of the
sequence $f_m,$ $m=1,2,\ldots .$

\medskip
Since $(\overline{D}, |\cdot|)$ is separable metric space and
$(\overline{D_0}_P, \rho)$ is a compact metric space, the normality
of the sequence $f_m,$ $m=1,2,\ldots ,$ follows by Arzela-Ascoli
theorem (see e.g. \cite[Theorem~20.4]{Va}). Thus,
$f_{m_k}\rightarrow f$ uniformly in $\overline{D}$ by the
metric~$\rho$ for some a subsequence of numbers $m_k,$ $k=1,2,\ldots
.$

\medskip
Observe that, $f$ has a continuous extension to $\overline{D}.$
Otherwise, due to the compactness of $\overline{D_0}_P$ there are
$x_0\in \partial D$ and at least two subsequences $x_m,
x^{\,\prime}_m\in D,$ $x_m, x^{\,\prime}_m\rightarrow x_0$ as
$m\rightarrow\infty,$ and $P_1, P_2\in \overline{D_0}_P,$ $P_1\ne
P_2,$ such that $f(x_m)\rightarrow P_1$ and
$f(x^{\,\prime}_m)\rightarrow P_2$ as $m\rightarrow \infty$ by the
metric $\rho.$ Now, $\rho(f(x_m), f(x^{\,\prime}_m))\geqslant
\varepsilon_0$ for some $\varepsilon_0>0$ and sufficiently large
$m.$ However, by the triangle inequality,
$$\rho(f(x_m), f(x^{\,\prime}_m))\leqslant $$$$\leqslant\rho(f(x_m), f_{m_k}(x_m))
+ \rho(f_{m_k}(x_m), f_{m_k}(x^{\,\prime}_m))+
\rho(f_{m_k}(x^{\,\prime}_m), f(x^{\,\prime}_m))\rightarrow 0$$
as $k\rightarrow \infty$ uniformly by $m,$ because $f_{m_k}$ is
equicontinuous family by the proving above and $f_{m_k}$ converges
to $f$ uniformly in $\overline{D}.$ Thus, $f$ has a continuous
extension in $\overline{D},$ as required.

\medskip
It remains to prove that, for any $x_0\in
\partial D$ there is $P_0\in E_{D_0}:=\overline{D_0}_P\setminus D_0$ such that, for any
$\varepsilon>0$ there is $\delta=\delta(\varepsilon)>0$ and
$M=M(\varepsilon)\in {\Bbb N}$ such that $\rho(f_{m_k}(x),
P_0)<\varepsilon$ for all $x\in B(x_0, \delta)\cap D$ and
$k\geqslant M_0.$ Now, we set $P_0=f(x_0).$ By the triangle
inequality, given $\varepsilon>0$ we obtain that
$$\rho(P_0, f_{m_k}(x))=$$
\begin{equation}\label{eq3A}
=\rho(f(x_0), f_{m_k}(x))\leqslant \rho(f(x_0), f_{m_k}(x_0))+
\rho(f_{m_k}(x_0), f_{m_k}(x))<\varepsilon\end{equation}
whenever $k\geqslant K=K(\varepsilon)$ and
$|x-x_0|<\delta=\delta(\varepsilon).$

\medskip
Observe that $P_0\in E_{D_0}:=\overline{D_0}_P\setminus D_0.$
Indeed, in the contrary case, when $P_0\in D_0,$ it follows
by~(\ref{eq3A}) that $f_{m_k}(x_{m_k})\rightarrow P_0$ as
$k\rightarrow\infty.$ Now, the sequence $y_k:=f_{m_k}(x_{m_k})$ lies
in some a compact $K$ in $D_0.$ However, by Lemma~\ref{lem1}
$h(x_{m_k},
\partial D)\geqslant h(f_{m_k}^{-1}(K), \partial D)\geqslant\delta>0$ for some
$\delta>0,$ that contradicts the definition of the sequence
$x_{m_k}.$ The obtain contradiction finishes the proof of
Theorem~\ref{th2}.~$\Box$

\medskip
{\it Proof of Corollary~\ref{cor1}.} Let $0<r_0:=\sup\limits_{y\in
D}|y-y_0|.$ We may assume that $Q$ is extended by zero outside
$D_0.$ Let $Q\in L^1(D_0).$ By the Fubini theorem (see, e.g.,
\cite[Theorem~8.1.III]{Sa}) we obtain that
$$\int\limits_{r_1<|y-y_0|<r_2}Q(y)\,dm(y)=\int\limits_{r_1}^{r_2}
\int\limits_{S(y_0, r)}Q(y)\,d\mathcal{H}^{n-1}(y)dr<\infty\,.$$
This means the fulfillment of the condition of the integrability of
the function $Q$ on the spheres with respect to any subset $E_1$ in
$[r_1, r_2].$~$\Box$

\medskip
In the context of Theorem~\ref{th2}, one might think that the family
of mappings~$f_m,$ $m=1,2,\ldots,$ always associates boundary points
$x_0$ of the domain $D$ with boundary points of the kernel $D_0,$
and inner points of $D$ with inner points of $D_0,$ but this is not
quite so. Let us consider the following instructive example.

\medskip
\begin{example}\label{ex1}
Following~\cite[Example~5.3]{SSD}, cf.~\cite[Proposition~6.3]{MRSY},
we put $p\geqslant 1$ such that $n/p(n-1)<1.$ Put $\alpha\in (0,
n/p(n-1)).$ Let $f_m: B(0, 2)\setminus\{0\}\rightarrow {\Bbb B}^n,$
$D:=B(0, 2)\setminus\{0\},$
$$f_m(x)\,=\,\left
\{\begin{array}{rr} \frac{(|x|-1)^{1/\alpha}}{|x|}\cdot x\,,
& 1+1/(m^{\alpha})\leqslant|x|\leqslant 2, \\
\frac{1/m}{1+(1/m)^{\alpha}}\cdot x\,, & 0<|x|<1+1/(m^{\alpha}) \ .
\end{array}\right.
$$
Using the approach applied in~\cite[Proposition~6.3]{MRSY} and
applying \cite[Theorems~8.1, 8.5]{MRSY} we may show that $f_m$
satisfy the relations~(\ref{eq2*A})--(\ref{eqA2}) in
$\overline{{\Bbb B}^n}$ with
$Q(x)=C\cdot\frac{1}{|x|^{\alpha(n-1)}}$ with some a constant $C>0$
whenever $Q\in L^1({\Bbb B}^n).$ Observe also that $f_m(B(0,
2)\setminus \{0\})={\Bbb B}^n\setminus \{0\},$ so ${\Bbb
B}^n\setminus \{0\}$ is a kernel of the sequence of domains
$f_m(B(0, 2)\setminus \{0\})$ trivially. The relation
$h(f_m^{-1}(E),
\partial D)\geqslant \delta>0$ obviously holds for infinitely many continua $E$ in $D_0.$

\medskip
Let now $x_0\in {\Bbb S}^{n-1}=S(0, 1)=\partial {\Bbb B}^n.$
Although $x_0\in {\rm Int\,}\bigl(B(0, 2)\setminus \{0\}\bigr),$
observe that $f_m(x_0)\rightarrow 0$ as $m\rightarrow \infty$ and
$0\in
\partial ({\Bbb B}^n\setminus\{0\}).$

\medskip
In the above example, the domain $D_0$ is not regular. In order for
the example to exactly correspond to the formulation of
Theorem~\ref{th2}, we may set $D=B^{+}(0, 2)=\{x\in B(0, 2):
x=(x_1,x_2,\ldots x_n), x_n>0\}$ and $D_0={\Bbb B}^{n+}=\{x\in {\Bbb
B}^n: x=(x_1,x_2,\ldots x_n), x_n>0\}.$ Now, $f_m$ mentioned above
map $B^{+}(0, 2)$ onto ${\Bbb B}^{n+},$ ${\Bbb B}^{n+}$ is a kernel
of a constant sequence of domains $f_m(B^{+}(0, 2))={\Bbb B}^{n+},$
$m=1,2,\ldots,$ and  $f_m$ satisfy the
relations~(\ref{eq2*A})--(\ref{eqA2}) in $\overline{{\Bbb B}^{n+}}$
with $Q(x)=C\cdot\frac{1}{|x|^{\alpha(n-1)}},$ $Q\in L^1({\Bbb
B}^{n+}).$ The relation $h(f_m^{-1}(E),
\partial D)\geqslant \delta>0$ holds for infinitely many continua $E$ in $D_0,$
as well. Obviously, the domain ${\Bbb B}^{n+}$ has a locally
quasiconformal boundary and, consequently, is a regular domain.
Observe that, the ball $B(0, 2)$ has a weakly flat boundary, see
e.g.~\cite[Theorem~4.5]{Na$_1$}. Let now $x_0\in {\Bbb S}^{n-1}\cap
{\rm Int}\,B^{+}(0, 2).$ Although $x_0\in {\rm Int}\,B^{+}(0, 2),$
observe that $f_m(x_0)\rightarrow 0$ as $m\rightarrow \infty$ and
$0\in
\partial {\Bbb B}^{n+}.$ All of the conditions of Theorem~\ref{th2}
are satisfied, in particular, $f_m\in \frak{F}^{E}_{Q, \delta}(D,
D_0)$ for $D=B(0, 2),$ $D_0={\Bbb B}^{n+},$
$Q=C\cdot\frac{1}{|x|^{\alpha(n-1)}}\in L^1({\Bbb B}^n)$ and some a
compactum $E\subset {\Bbb B}^{n+}$ and $\delta>0.$ By this theorem,
$f_{m_k}(x)\rightarrow P_0$ as $k\rightarrow\infty$ for some a
subsequence of numbers $m_k,$ $k=1,2,\ldots $ (as we see, not only
$f_{m_k},$ but itself a sequence $f_m$ has the above property,
because the Euclidean convergence $f_m(x_0)\rightarrow 0$ as
$m\rightarrow \infty$ implies trivially the convergence by the
metric $\rho$ in any domain with a locally quasiconformal boundary).
By the reasons mentioned above, $P_0\in E_{{\Bbb B}^{n+}}.$
\end{example}

\medskip
Based on Example~\ref{ex1}, an inner point in the preimage $D$ may
correspond to either a prime end of the kernel $D_0$ or an inner
point of this kernel. At the same time, according to
Theorem~\ref{th2}, a prime end in the preimage can only correspond
to a prime end of the kernel. Let us ask the problem: when does this
correspondence look more ``equal''? The answer to this question is
given below.

\medskip
\begin{theorem}\label{th1}
{\it\, Assume that, under assumptions and notions of
Theorem~\ref{th2}, $f_m$ and $Q$ satisfy some additional conditions:
$h(f_m(E_2))\geqslant \delta_2$ for some a continuum $E_2$ in $D,$
some $\delta_2>0$ and every $m\in {\Bbb N}.$ Besides that, let $Q\in
L^1(D_0)$ and at least one of the following conditions hold:

\medskip
1) $Q\in FMO(\overline{D_0});$

\medskip
2) for any $y_0\in \overline{D_0}$ there is $\delta(y_0)>0$ such
that, $\int\limits_{\varepsilon}^{\delta(y_0)}
\frac{dt}{tq_{y_0}^{\frac{1}{n-1}}(t)}<\infty$ for any
$\varepsilon\in (0, \delta(y_0))$ and
\begin{equation}\label{eq5D}
\int\limits_{0}^{\delta(y_0)}
\frac{dt}{tq_{y_0}^{\frac{1}{n-1}}(t)}=\infty\,.
\end{equation}
Then under notions of the statement of Theorem~\ref{th2}, the
following condition is true: for any $x_0\in D$ there is
$y_0:=f(x_0)\in D_0$ such that, for any $\varepsilon>0$ there is
$\delta=\delta(\varepsilon)>0$ and $M=M(\varepsilon)\in {\Bbb N}$
such that $\rho(f_{m_k}(x), y_0)<\varepsilon$ for all $x\in B(x_0,
\delta)\cap D$ and $k\geqslant M_0.$ Moreover, $f(D)=D_0$ and
$f(\overline{D})=\overline{D_0}_P.$ }
\end{theorem}

\medskip
\begin{proof}
{\bf I.} Let $f_{m_k},$ $k=1,2,\ldots,$ be a sequence which
converges to $f$ uniformly by the metric $\rho$ as
$k\rightarrow\infty.$ Let $x_0\in D.$ Let us to prove the first part
of the statement of the theorem: for any $x_0\in D$ there is
$y_0=f(x_0)\in D_0$ such that, for any $\varepsilon>0$ there is
$\delta=\delta(\varepsilon)>0$ and $M=M(\varepsilon)\in {\Bbb N}$
such that $\rho(f_{m_k}(x), y_0)<\varepsilon$ for all $x\in B(x_0,
\delta)\cap D$ and $k\geqslant M_0.$ If $y_0=f(x_0)\in D_0,$ there
is nothing to prove. Let us show that the case $y_0\in E_{D_0}$ is
impossible. We prove the latter by the contradiction, i.e., assume
that $y_0\in E_{D_0}.$ Since $\overline{D_0}$ is a compactum in
${\Bbb R}^n,$ there is a subsequence $f_{{m_k}_l}(x_0),$
$l=1,2,\ldots,$ converging to some point $z_0\in \overline{D_0}.$ In
order not to complicate the notation, we will assume that the
sequence $f_{m_k}(x_0),$ $k=1,2,\ldots,$ itself has this property,
i.e., $f_{m_k}(x_0)\rightarrow z_0$ as $k\rightarrow\infty$ by the
Euclidean metric. Now, $z_0\in \partial D_0$ because $f_{m_k}(x_0)$
converges to $y_0\in E_{D_0}$ as $k\rightarrow\infty$ by the
assumption.

\medskip
On the other hand, let $B(x_0, \varepsilon_1)\subset D$ for some
$\varepsilon_1>0,$ moreover, $\overline{B(x_0,
\varepsilon_1)}\subset D.$ By Proposition~\ref{pr3}, there is
$r_0>0,$ which does not depend on $m\in {\Bbb N},$ such that
\begin{equation}\label{eq1B}
f_{m_k}(B(x_0, \varepsilon_1))\supset B(f_{m_k}(x_0), r_0)\qquad
\forall\,\, k\in{\Bbb N}\,.
\end{equation}
Here we took into account that, obviously, $B(f_m(x_0), r_0)\subset
B_h(f_m(x_0), r_0)=\{y\in \overline{{\Bbb R}^n}: h(f_m(x_0),
y)<r_0\}.$ Since $f_{m_k}(x_0)\rightarrow z_0$ as
$k\rightarrow\infty,$ we have that $z_0\in B(f_{m_k}(x_0), r_0)$ for
sufficiently large $k\in {\Bbb N}.$ By~(\ref{eq1B}), $z_0\in
f_{m_k}(B(x_0, \varepsilon_1))\subset D_0.$ This contradicts the
assumption $z_0\in\partial D_0$ made above. The obtained
contradiction completes the first part of the proof of
Theorem~\ref{th1}.

\medskip
{\bf II.} Let us to prove that $f(D)=D_0.$ We may consider that
$f_m$ converges to $f$ uniformly in $D$ as $m\rightarrow\infty.$ If
$x_0\in D,$ then $f(x_0)=\lim\limits_{m\rightarrow\infty}f_m(x_0),$
where ``{\rm lim}'' must be understood in the metrics $\rho.$ By the
proved above, $f(x_0)\in D_0,$ thus $f(D)\subset D_0.$ On the other
hand, let $y_0\in D_0.$ Now, since $D_0$ is a kernel of $f_m(D),$ we
have that $y_0=f_m(x_m)$ for sufficiently large $m=1,2,\ldots .$ We
may consider that $x_m\rightarrow x_0\in\overline{D}$ as
$m\rightarrow\infty.$ By the triangle inequality
$$\rho(f_m(x_m), f(x_0))\leqslant
\rho(f_m(x_m), f(x_m))+\rho(f(x_m), f(x_0))\,.$$
Since $f_m$ converges to $f$ uniformly in $D$ and $f$ has a
continuous extension to $x_0$ (see Theorem~\ref{th2}), the latter
implies that $y_0=f_m(x_m)\rightarrow f(x_0).$ Thus, $y_0=f(x_0).$
Since $y_0\in D_0,$ we have that $x_0\in D$ (see Theorem~\ref{th2}).
Now, $y_0\in f(D),$ so that $D_0\subset f(D).$ Thus, $f(D)=D_0,$ as
required.

\medskip
{\bf III.} Finally, let us to prove that
$f(\overline{D})=\overline{D_0}_P.$ If $x_0\in \overline{D},$ then
either $x_0\in D,$ or $x_0\in\partial D.$ If $x_0\in D,$ then
$f(x_0)\in D_0\subset \overline{D_0}_P$ by the proven above. If
$x_0\in \partial D,$ then $f(x_0)\in E_{D_0}\subset\overline{D_0}_P$
by Theorem~\ref{th2}. In any of two cases, $f(x_0)\in
\overline{D_0}_P,$ therefore,
$f(\overline{D})\subset\overline{D_0}_P.$ Otherwise, let $y_0\in
\overline{D_0}_P.$ There are two cases: $y_0\in D_0$ or $y_0\in
E_{D_0}.$ In the first case, when $y_0\in D_0,$ we have that
$y_0=f(x_0)$ for some point $x_0\in D$ (see the step~{\bf{II}}) and,
consequently, $y_0\in f(\overline{D}).$ In the second case, when
$y_0\in E_{D_0},$ there is a sequence $y_k\in D_0,$ $y_k\rightarrow
y_0$ as $k\rightarrow\infty$ in the metric $\rho.$ By the definition
of the kernel $D_0,$ given $k\in {\Bbb N}$ there is $m_k$ such that
$y_k\in f_m(D)$ for any $m\geqslant m_k.$ We may consider that the
sequence $m_k$ is increasing by $k.$ Now, $y_k=f_{m_k}(x_k)$ for
some $x_k\in D.$ We may consider that $x_k\rightarrow x_0$ as
$k\rightarrow \infty,$ where $x_0\in \overline{D}.$ Since $f_m$
converges uniformly to $f$ as $m\rightarrow\infty,$ we obtain that
$$\rho(f_{m_k}(x_k), f(x_0))\leqslant
\rho(f_{m_k}(x_k), f(x_k))+\rho(f(x_k), f(x_0))\rightarrow 0$$
as $k\rightarrow\infty.$ It follows from that,
$y_k=f_{m_k}(x_k)\rightarrow f(x_0)$ and simultaneously
$y_k\rightarrow y_0$ as $k\rightarrow\infty.$ Thus $y_0=f(x_0)\in
f(\overline{D}).$ So, we have proved that $\overline{D_0}_P\subset
f(\overline{D}).$ This proved the equality
$f(\overline{D})=\overline{D_0}_P.$ Theorem is proved.~$\Box$
\end{proof}

\section{Lemma on approaching continua}

Recall that $\frak{F}^{E}_{Q, \delta}(D, D_0)$ denotes the class of
all open, discrete and closed mappings $f:D\rightarrow D_0$ of a
domain $D$ onto some domain $f(D),$ $E\subset f(D)\subset D_0,$ such
that (\ref{eq2*A})--(\ref{eqA2}) hold for every $y_0\in
\overline{D_0}$ and, in addition, $f^{\,-1}(E)$ is a continuum such
that $h(f^{\,-1}(E),
\partial D)\geqslant \delta.$ Besides that, $\frak{R}^{E}_{Q, \delta}(D, D_0)$
is a class of all open, discrete and closed mappings $f:D\rightarrow
D_0$ of a domain $D$ onto some domain $f(D),$ $E\subset f(D)\subset
D_0,$ such that (\ref{eq2*A})--(\ref{eqA2}) hold for every $y_0\in
\overline{D_0}$ and, in addition, $f^{\,-1}(E)$ is a continuum with
$h(f^{\,-1}(E))\geqslant \delta.$ The following statement holds,
cf.~\cite[Lemma~4.1]{SevSkv$_2$}, \cite[Lemma~2.13]{ISS}.

\medskip
\begin{lemma}\label{lem2}
{\it\, Let $\delta>0,$ let $D, D_0$ be domains in ${\Bbb R}^n,$
$n\geqslant 2,$ let $E$ be a continuum in $D_0$ and let $Q:{\Bbb
R}^n\rightarrow [0, \infty]$ be a Lebesgue measurable function,
$Q(y)\equiv 0$ outside $D_0.$ Assume that, for each point $y_0\in
\overline{D_0}$ and for every $0<r_1<r_2<r_0:=\sup\limits_{y\in
D_0}|y-y_0|$ there is a set $E_1\subset[r_1, r_2]$ of a positive
linear Lebesgue measure such that the function $Q$ is integrable
with respect to $\mathcal{H}^{n-1}$ over the spheres $S(y_0, r)$ for
every $r\in E_1.$ Assume that, no connected component of the
boundary of the domain $D$ degenerates into a point and, besides
that, $f_{m}(D)$ converge to $D_0$ as its kernel as
$m\rightarrow\infty$ for $f_m\in \frak{R}^{E}_{Q, \delta}(D, D_0),$
$m=1,2,\ldots .$

Now, there is $\delta_*>0$ such that $h(f_m^{\,-1}(E),
\partial D)\geqslant \delta_*$ for every $m\in {\Bbb N},$ i.e.,
$f_m\in \frak{F}^{E}_{Q, \delta_*}(D, D_0).$ }
\end{lemma}

\medskip
\begin{proof}
The proof of Lemma~\ref{lem2} is in many ways similar to
Lemma~\ref{lem1}. Note also that, under open, discrete, and closed
mappings, the preimage of any compact set $E$ is compact
(see~\cite[Theorem~3.3]{Vu}), so that the presence of the number
$\delta_*>0$ in condition $h(f_m^{\,-1}(E),
\partial D)\geqslant \delta_*$ is obvious for any fixed $m\in {\Bbb N}$
and some $\delta_*=\delta_*(m).$ The question is only about the
presence of a common delta $\delta_*>0$ that provides the entire
family of mappings $\frak{R}^{E}_{Q, \delta}(D, D_0).$

\medskip
Let us prove Lemma~\ref{lem2} by the contradiction. Assume that, the
conclusion of the lemma is not true. Then for each $k\in {\Bbb N}$
there is some number $m_k\in {\Bbb N}$ such that
$h(f^{\,-1}_{m_k}(E),
\partial D)<1/m_k.$ We may assume that the sequence $m_k$ is
increasing by $k\in {\Bbb N},$ moreover, we may consider that the
latter holds for any $m\in {\Bbb N},$ i.e., for each $m\in {\Bbb N}$
there is some number $m\in {\Bbb N}$ such that $h(f^{\,-1}_{m}(E),
\partial D)<1/m.$ Since $f^{\,-1}_{m}(E)$ is a compactum in $D$ for each
$m=1,2,\ldots,$ we have that $h(f^{\,-1}_{m}(E),
\partial D)=h(x_m, y_m)<\frac{1}{m}$ for some
$x_m\in f^{\,-1}_{m}(E)$ and $y_m\in \partial D.$ Since $\partial D$
is a compact set, we may assume that $y_m\rightarrow y_0$ as
$m\rightarrow\infty$ for some $y_0\in
\partial D.$ Now, $x_m\rightarrow y_0$ as $m\rightarrow\infty,$ as
well.

\medskip
Set $A_m:=f^{\,-1}_{m}(E).$ Let $E_0$ be a component of $\partial D$
consisting $y_0.$ Since by the assumptions of the lemma, all
components of $\partial D$ are non-degenerate, there exists $r>0$
such that $h(E_0)\geqslant r.$ Put $P>0$ and $U=B_h(y_0, R_0)=\{y\in
\overline{{\Bbb R}^n}: h(y, y_0)<R_0\},$ where $2R_0:=\min\{r/2,
\delta/2\}$ and $\delta$ is a number in the definition of the class
$\frak{R}^{E}_{Q, \delta}(D, D_0)$. Observe that $A_m\cap
U\ne\varnothing\ne A_m\setminus U$ for sufficiently large $m\in{\Bbb
N},$ since $x_m\rightarrow y_0$ as $m\rightarrow \infty,$ $x_m\in
A_m;$ besides that, $h(A_m)\geqslant \delta\geqslant 2R_0$ and
$h(U)\leqslant 2R_0.$ Since $A_m$ is a continuum, $A_m\cap
\partial U\ne\varnothing$ by Proposition~\ref{pr2}. Similarly,
$E_0\cap U\ne\varnothing\ne E_0\setminus U$ for sufficiently large
$m\in{\Bbb N},$ since $h(E_0)\geqslant r> 2R_0$ and $h(U)\leqslant
2R_0.$ Since $E_0$ is a continuum, $E_0\cap
\partial U\ne\varnothing$ by Proposition~\ref{pr2}.
By the proving above,
\begin{equation}\label{eq8A}
A_m\cap
\partial U\ne\varnothing\ne E_0\cap
\partial U\,.
\end{equation}
By Lemma~\ref{lem3} there is $V\subset U,$ $V$ is a neighborhood of
$y_0,$ such that
\begin{equation}\label{eq9A}
M(\Gamma(E, F, \overline{{\Bbb R}^n}))>P
\end{equation}
for any continua $E, F\subset \overline{{\Bbb R}^n}$ with $E\cap
\partial U\ne\varnothing\ne E\cap \partial V$ and $F\cap \partial
U\ne\varnothing\ne F\cap \partial V.$
Arguing similarly to above, we may prove that
$$
A_m\cap
\partial V\ne\varnothing\ne E_0\cap
\partial V
$$
for sufficiently large $m\in {\Bbb N}.$ Thus, by~(\ref{eq9A})
\begin{equation}\label{eq9D}
M(\Gamma(A_m, E_0,\overline{{\Bbb R}^n}))>P
\end{equation}
for sufficiently large $m=1,2,\ldots .$ Let $\gamma:[0,
1]\rightarrow \overline{{\Bbb R}^n}$ be a path in $\Gamma(A_m, E_0,
\overline{{\Bbb R}^n}),$ i.e., $\gamma(0)\in A_m,$ $\gamma(1)\in
E_0$ and $\gamma(t)\in \overline{{\Bbb R}^n}$ for $t\in (0, 1).$ Let
$t_m=\sup\limits_{\gamma(t)\in D}t$ and let
$\alpha_m(t)=\gamma|_{[0, t_m)}(t).$ Let $\Gamma_m$ consists of all
such paths $\alpha_m,$ now $\Gamma(A_m, E_0, \overline{{\Bbb
R}^n})>\Gamma_m$ and by the minorization principle of the modulus
(see~\cite[Theorem~1]{Fu})
\begin{equation}\label{eq12A}
M(\Gamma_m)\geqslant M(\Gamma(A_m, E_0, \overline{{\Bbb R}^n}))\,.
\end{equation}
Combining~(\ref{eq9D}) and~(\ref{eq12A}), we obtain that
\begin{equation}\label{eq9E}
M(\Gamma_m)>P
\end{equation}
for sufficiently large $m=1,2,\ldots .$

\medskip
We now prove that the relation~(\ref{eq9E}) contradicts the
definition of $f_m$ in~(\ref{eq2*A})--(\ref{eqA2}). Since $f_m$ is a
closed mapping, it preserves the boundary (see
\cite[Theorem~3.3]{Vu}), so that $C(f_m,
\partial D)\subset \partial f_m(D).$ Thus, $C(\beta_m(t), t_m)\subset \partial f_m(D),$
$\beta_m(t):=f_m(\alpha_m(t))=f_m(\gamma|_{[0, t_m)}(t)),$
$\gamma\in \Gamma(A_m, E_0, \overline{{\Bbb R}^n})$ and
$t_m=\sup\limits_{\gamma(t)\in D}t.$

\medskip
We cover the continuum $E$ with balls $B(x, \varepsilon/4),$ $x\in
A.$ Since $E$ is a compact set, we may assume that $E\subset
\bigcup\limits_{i=1}^{M_0}B(x_i, \varepsilon/4),$ $x_i\in K_0,$
$i=1,2,\ldots, M_0,$ $1\leqslant M_0<\infty.$ By the definition,
$M_0$ depends only on $E,$ in particular, $M_0$ does not depend on
$m.$
Note that
\begin{equation}\label{eq6E}
\Gamma_m=\bigcup\limits_{i=1}^{M_0}\Gamma_{mi}\,,
\end{equation}
where $\Gamma_{mi}$ consists of all paths $\gamma:[0, 1)\rightarrow
D$ in $\Gamma_m$ such that $f_m(\gamma(0))\in B(y_i,
\varepsilon/4).$ As under the proof of Lemma~\ref{lem1}, we may show
that
\begin{equation}\label{eq7I}
f_m(\Gamma_{mi})>\Gamma(S(y_i, \varepsilon/4), S(y_i,
\varepsilon/2), A(y_i, \varepsilon/4, \varepsilon/2))\,.
\end{equation}
By~(\ref{eq6E}) and (\ref{eq7I}), we obtain that
\begin{equation}\label{eq5A}
\Gamma_m>\bigcup\limits_{i=1}^{M_0}\Gamma_{f_m}(y_i, \varepsilon/4,
\varepsilon/2)\,.
\end{equation}
Set $\widetilde{Q}(y)=\max\{Q(y), 1\}$ and
$$\widetilde{q}_{y_i}(r)=\int\limits_{S(y_i,
r)}\widetilde{Q}(y)\,d\mathcal{A}\,.$$ Now,
$\widetilde{q}_{y_i}(r)\ne \infty$ for $r\in
E_1\subset[\varepsilon/4, \varepsilon/2],$ where $E_1$ is some set
of positive linear measure which exists by the assumptions of the
lemma. Set
$$I_i=I_i(y_i, \varepsilon/4, \varepsilon/2)=\int\limits_{\varepsilon/4}^{\varepsilon/2}\
\frac{dr}{r\widetilde{q}_{y_i}^{\frac{1}{n-1}}(r)}\,.$$
Observe that $I\ne 0,$ because $\widetilde{q}_{y_i}(r)\ne \infty$
for $r\in E_1\subset[\varepsilon/4, \varepsilon/2],$ where $E_1$ is
some set of positive linear measure. Besides that, by the direct
calculations we obtain that
$$I_i\leqslant \log\frac{r_2}{r_1}<\infty\,,\quad i=1,2, \ldots, M_0\,.$$
Now, we put
$$\eta_i(r)=\begin{cases}
\frac{1}{I_ir\widetilde{q}_{y_i}^{\frac{1}{n-1}}(r)}\,,&
r\in [\varepsilon/4, \varepsilon/2]\,,\\
0,& r\not\in [\varepsilon/4, \varepsilon/2]\,.
\end{cases}$$
Observe that, a function~$\eta_i$ satisfies the
condition~$\int\limits_{\varepsilon/4}^{\varepsilon/2}\eta_i(r)\,dr=1,$
therefore it may be substituted into the right side of the
inequality~(\ref{eq2*A}) with the corresponding values $f,$ $r_1$
and $r_2.$ Now, we obtain that
\begin{equation}\label{eq7J}
M(\Gamma_{f_m}(y_i, \varepsilon/4, \varepsilon/2))\leqslant
\int\limits_{A(y_i, \varepsilon/4, \varepsilon/2)}
\widetilde{Q}(y)\,\eta^n_i(|y- y_i|)\,dm(y)\,.\end{equation}
By the Fubini theorem we have that
$$\int\limits_{A(y_i, \varepsilon/4, \varepsilon/2)}
\widetilde{Q}(y)\,\eta^n_i(|y- y_i|)\,dm(y)=$$
\begin{equation}\label{eq7K}
=\int\limits_{\varepsilon/4}^{\varepsilon/2}\int\limits_{S(y_i,
r)}Q(y)\eta^n_i(|y- y_i|)\,d\mathcal{A}\,dr\,=
\end{equation}
$$=\frac{\omega_{n-1}}{I_i^n}\int\limits_{\varepsilon/4}^{\varepsilon/2}r^{n-1}
\widetilde{q}_{y_i}(r)\cdot
\frac{dr}{r^n\widetilde{q}^{\frac{n}{n-1}}_{y_i}(r)}=\frac{\omega_{n-1}}{I_i^{n-1}}\,,$$
where $\omega_{n-1}$ is the area of the unit sphere in ${\Bbb R}^n.$
Now, by~(\ref{eq7J}) and~(\ref{eq7K}) we obtain that
$$M(\Gamma_{f_m}(y_i, \varepsilon/4,
\varepsilon/2))\leqslant \frac{\omega_{n-1}}{I_i^{n-1}}\,,$$
whence from~(\ref{eq5A}) we obtain that
\begin{equation}\label{eq7L}
M(\Gamma_m)\leqslant \sum\limits_{i=1}^{M_0}M(\Gamma_{f_m}(y_i,
\varepsilon/4, \varepsilon/2))\leqslant
\sum\limits_{i=1}^{M_0}\frac{\omega_{n-1}}{I_i^{n-1}}:=C_0\,, \quad
m=1,2,\ldots\,.
\end{equation}
Since $P$ in~(\ref{eq9E}) may be done arbitrary big, the
relations~(\ref{eq9E}) and~(\ref{eq7L}) contradict each other. This
completes the proof.~$\Box$
\end{proof}

\section{Proof of another results and some examples}

{\it Proof of Theorem~\ref{th3}} directly follows by
Theorem~\ref{th2} and Lemma~\ref{lem2}.~$\Box$

\medskip
{\it Proof of Corollary~\ref{cor2}} directly follows by
Theorem~\ref{th3} on the basis of arguments given under the proof of
Corollary~\ref{cor1}.~$\Box$

\medskip
Similarly to Theorem~\ref{th1} we may formulate the following
statement.

\medskip
\begin{theorem}\label{th4}
{\it\, Assume that, under assumptions and notions of
Theorem~\ref{th3}, $f_m$ and $Q$ satisfy some additional conditions:
$h(f_m(E_2))\geqslant \delta_2$ for some a continuum $E_2$ in $D,$
some $\delta_2>0$ and every $m\in {\Bbb N}.$ Besides that, let $Q\in
L^1(D_0)$ and at least one of the following conditions hold:

\medskip
1) $Q\in FMO(\overline{D_0});$

\medskip
2) for any $y_0\in \overline{D_0}$ there is $\delta(y_0)>0$ such
that, $\int\limits_{\varepsilon}^{\delta(y_0)}
\frac{dt}{tq_{y_0}^{\frac{1}{n-1}}(t)}<\infty$ for any
$\varepsilon\in (0, \delta(y_0))$ and
$$
\int\limits_{0}^{\delta(y_0)}
\frac{dt}{tq_{y_0}^{\frac{1}{n-1}}(t)}=\infty\,.
$$
Then under notions of the statement of Theorem~\ref{th3}, the
following condition is true: for any $x_0\in D$ there is
$y_0:=f(x_0)\in D_0$ such that, for any $\varepsilon>0$ there is
$\delta=\delta(\varepsilon)>0$ and $M=M(\varepsilon)\in {\Bbb N}$
such that $\rho(f_{m_k}(x), y_0)<\varepsilon$ for all $x\in B(x_0,
\delta)\cap D$ and $k\geqslant M_0.$ Moreover, $f(D)=D_0$ and
$f(\overline{D})=\overline{D_0}_P.$ }
\end{theorem}

\medskip
{\it Proof of Theorem~\ref{th4}} directly follows by
Theorem~\ref{th1} and Lemma~\ref{lem2}.~$\Box$

\medskip
\begin{example}\label{ex2}
Let $D_0$ be the unit disk ${\Bbb D}$ in ${\Bbb C}$ with a cut along
the segment $[0, 1]$, i.e., $D_0={\Bbb D}\setminus I \subset {\Bbb
C}$, where ${\Bbb D}:=\{z\in {\Bbb C}\,:\, |z|<1\}$, $I:=\{z=x+iy\in
{\Bbb C}\,:\, y=0, 0\leqslant  x\leqslant  1\}$. Note that $D_0$ is
a regular domain because by the Riemannian mapping theorem, $D_0$ is
conformally equivalent to ${\Bbb D}$; moreover, ${\Bbb D}$ has
locally quasiconformal boundary and, consequently, weekly flat
boundary (see e.g. \cite[Theorem~17.10]{Va}). Let $f_m$ be a
conformal mapping of ${\Bbb D}$ onto $D_m=\{z\in {\Bbb C}\,:\,
|z|<m/(m+1)\}\setminus I_m,$ $I_m:=\{z=x+iy\in {\Bbb C}\,:\, y=0,
0\leqslant  x\leqslant m/(m+1)\}$. The above mapping $f_m$ also
exists by Riemannian mapping theorem, moreover, we may consider that
$f_m(0)=-1/4.$

\medskip
Observe that $D_m$ converge to $D_0={\Bbb D}\setminus I$ as its
kernel, see Figure~\ref{fig2}.
\begin{figure}
  \centering\includegraphics[scale=0.5]{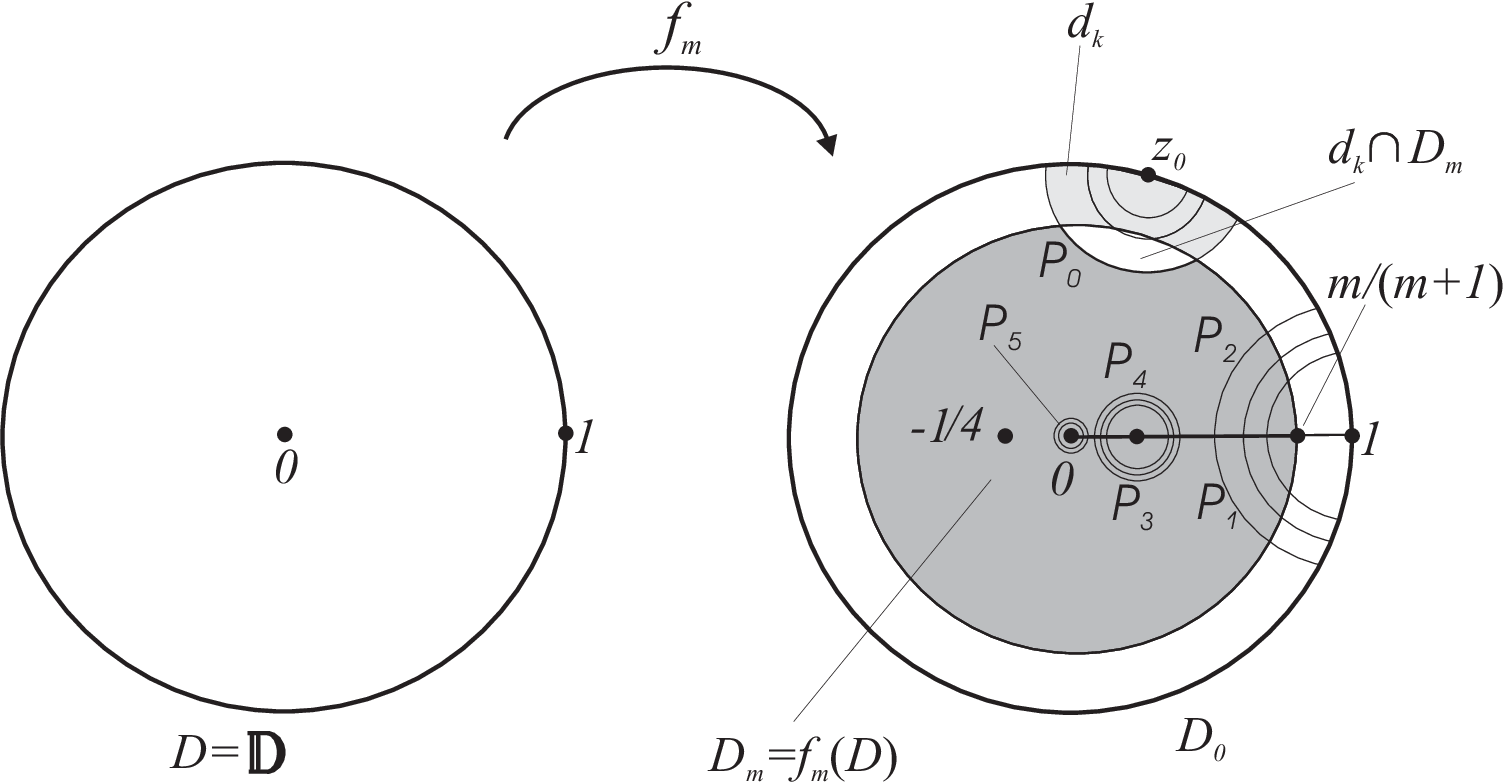}
  \caption{Illustration for Example~\ref{ex2}.}\label{fig2}
 \end{figure}
Besides that, the domains $d_k,$ formed by chains of cuts lying on
the circles centered at some point $z_0\in
\partial {\Bbb D}\cup I$ form a connected intersection with the
domains $D_m$ for each fixed $k$ and all sufficiently large $m.$
Indeed, prime ends of $D_0$ may be conditionally divided into~6
groups: 0) the set of prime ends $P_0$ which correspond to $z_0\in
{\Bbb S}^1\setminus I,$ where ${\Bbb S}^1=\partial {\Bbb D};$ 1) and
2): the prime ends $P_1$ and $P_2$ with impression $z_0=1,$ the
sequences of domains $d_k$ lie from below and from upper of $I,$
correspondingly; 3) and 4): the set of prime ends $P_3$ and $P_4$
with impression at some $z_0\in I\setminus \{1,0\},$ the sequences
of domains $d_k$ lie from below and from upper of $I,$
correspondingly; 5) the prime end $P_5$ which correspond to its
impression $I(P_5)=\{0\}.$ Let $d_k,$ $k=1,2,\ldots,$ be a sequence
of domains in $P_i,$ $i=\overline{0, 5}.$ Now, $d_k$ is convex for
any $k=1,2,\ldots. $ Given $k\in {\Bbb N},$ the intersection
$D_m\cap d_k$ is non-empty and convex for sufficiently large $m.$ In
particular, $D_m\cap d_k$ is connected for above $m,$ as required.

\medskip
All mappings $f_m$ are conformal and, consequently,
satisfy~(\ref{eq2*A})--(\ref{eqA2}) with $Q\equiv 1$ (see
\cite[Theorem~8.1]{Va}). In addition, $f_m$ has a continuous
extension to $\overline{D}$ by the metric $\rho,$ because $f_m$ has
a usual (Euclidean) continuous extension to $\overline{D}$ (see e.g.
\cite[Theorem~4.2]{Na$_1$}), in addition, $D_m$ are compact
subdomains of $D_0$ and metric $\rho$ in $\overline{D_0}_P$ is
homeomorphic to the Euclidean metric on compact sets of $D_0.$ Thus,
all of the conditions of Theorem~\ref{th2} are satisfied, in
particular, we may set $E=\{-1/4\}$ and $\delta:=h(0, \partial {\Bbb
D}).$ By the same reasons, all of the conditions of
Theorem~\ref{th1} hold for some continuum $E_2$ and some
$\delta_2>0,$ and conclusions of this theorem hold, as well.
\end{example}

\medskip
\begin{example}\label{ex3}
Similar example may be constructed in the space, as well. Let $x=(z,
x_{n-1}, x_n)\in {\Bbb R}^n$, where $z=(x_1,x_2,\dots , x_{n-2})\in
{\Bbb  R}^{n-2}$. Let $(r, \varphi)$ be the polar coordinates of the
point $(x_{n-1}, x_n)$: $x_{n-1}=r\cos\varphi$, $x_n=r\sin\varphi$,
$r=\sqrt{x_1^2+x_2^2}$. We set $D_0=\{x=(z, x_{n-1}, x_n)\in {\Bbb
B}^n\,:\, 0<\varphi<2\pi\},$ $D={\Bbb B}^n.$ We define $f_0(x)=(z,
r\cos\varphi/2, r\sin\varphi/2)$. Due to \cite[Example~16.3]{Va},
$f_0$ is a quasiconformal mapping of~$D_0$ onto the half-ball
$$
{{\Bbb  B}^{n+}=\{x=(z, x_1, x_2)\in {\Bbb  B}^n\,:\,
0<\varphi<\pi\}}\,,
$$
which has a locally quasiconformal boundary
(see~\cite[Lemma~40.2]{Va}). Thus, $D_0$ is a regular domain.
Observe that, ${\Bbb B}^n$ has a weakly flat boundary, see e.g.
\cite[Theorems~17.10, 17.12]{Va}.

\medskip
Below we construct some quasiconformal mapping of ${\Bbb B}^n$ onto
${\Bbb B}^{n+}.$ First, we define the mapping
$$f(x)=\left(x_1, x_2, \ldots, x_{n-1},
\frac{\sqrt{1-(x_1^2+\ldots+x_{n-1}^2)}+x_n}{2}\right)\,,$$
transforming ${\Bbb B}^n$ onto ${\Bbb B}^{n+}=\{x=(x_1,\ldots
x_n)\in {\Bbb B}^n: x_n>0\},$ whose Jacobian is 1/2. This mapping is
not quasiconformal, since the derivatives $\frac{\partial
f_n}{\partial x_i}=-\frac{x_i}{2\sqrt{1-(x_1^2+\ldots+x_{n-1}^2)}}$
are not bounded near of the ``equator''
$E:=\{x_1^2+\ldots+x_{n-1}^2=1, x_n=0\}$ of the sphere ${\Bbb
S}^{n-1}.$ In order to ``correct'' the mapping $f,$ we first cut off
the ball ${\Bbb B}^n$ so as to remove the ``bad'' set $E.$ Given
$0<h<1,$ consider a mapping
$$g(x)=\quad\left\{
\begin{array}{rr}
\left(\frac{h}{\sqrt{1-x_n^2}}x_1, \frac{h}{\sqrt{1-x_n^2}}x_2,
\ldots, \frac{h}{\sqrt{1-x_n^2}}x_{n-1}, x_n\right), & x^2_n< 1-h^2,\\
x, & x^2_n\geqslant 1-h^2\end{array} \right.\,.
$$
Observe that, $g$ maps ${\Bbb B}^n$ onto the ``cut ball''
$B_h=\{x\in {\Bbb B}^n: |x_i|<h, i=1,2,\ldots, n-1\}.$ We show that
$g$ is quasiconformal. Indeed, $g|_{S_y}$ is a radial mapping on a
fixed sphere $S_y=\{x_n=y\},$ in other words,
$g_{S_y}(\widetilde{x})=\frac{\widetilde{x}}{|\widetilde{x}|}\rho(|\widetilde{x}|),$
where $\widetilde{x}=(x_1,\ldots, x_{n-1})$. Using approaches and
terminology applied under the consideration of Proposition~6.3 in
\cite{MRSY}, we have that
$$\delta_{\tau}=\frac{|g_{S_y}(\widetilde{x})|}{|\widetilde{x}|}=\delta_r=\frac{\partial
|g_{S_y}|}{\partial(|\widetilde{x}|)}=\frac{h}{\sqrt{1-x_n^2}}\,.$$
Now, $\delta_{\tau}=\delta_r<1$ for $x^2_{n-1}<1-h^2,$ so that
$\Vert g^{\,\prime}(x)\Vert=1,$ $|J(x,
g)|=\left(\frac{h}{\sqrt{1-x_n^2}}\right)^{n-1}$ for
$x^2_{n-1}<1-h^2.$ If $x^2_{n-1}\geqslant 1-h^2,$ then $\Vert
g^{\,\prime}(x)\Vert=J(x, g)=1.$ Thus,
$$K_O(x, g)=\frac{\Vert g^{\,\prime}(x)\Vert^n}{|J(x, g)|}=
\left(\frac{\sqrt{1-x_n^2}}{h}\right)^{n-1}\leqslant\frac{1}{h^{n-1}}\,.$$
Finally, the mapping $\varphi(x)=g^{-1}\circ f\circ g$ transforms
${\Bbb B}^n$ onto ${\Bbb B}^n_+=\{x=(x_1,\ldots x_n)\in {\Bbb B}^n:
x_n>0\}$ and $\varphi$ is a quasiconformal because all the mappings
that make it up, are quasiconformal, see Figure~\ref{fig3}.
\begin{figure}
  \centering\includegraphics[scale=0.5]{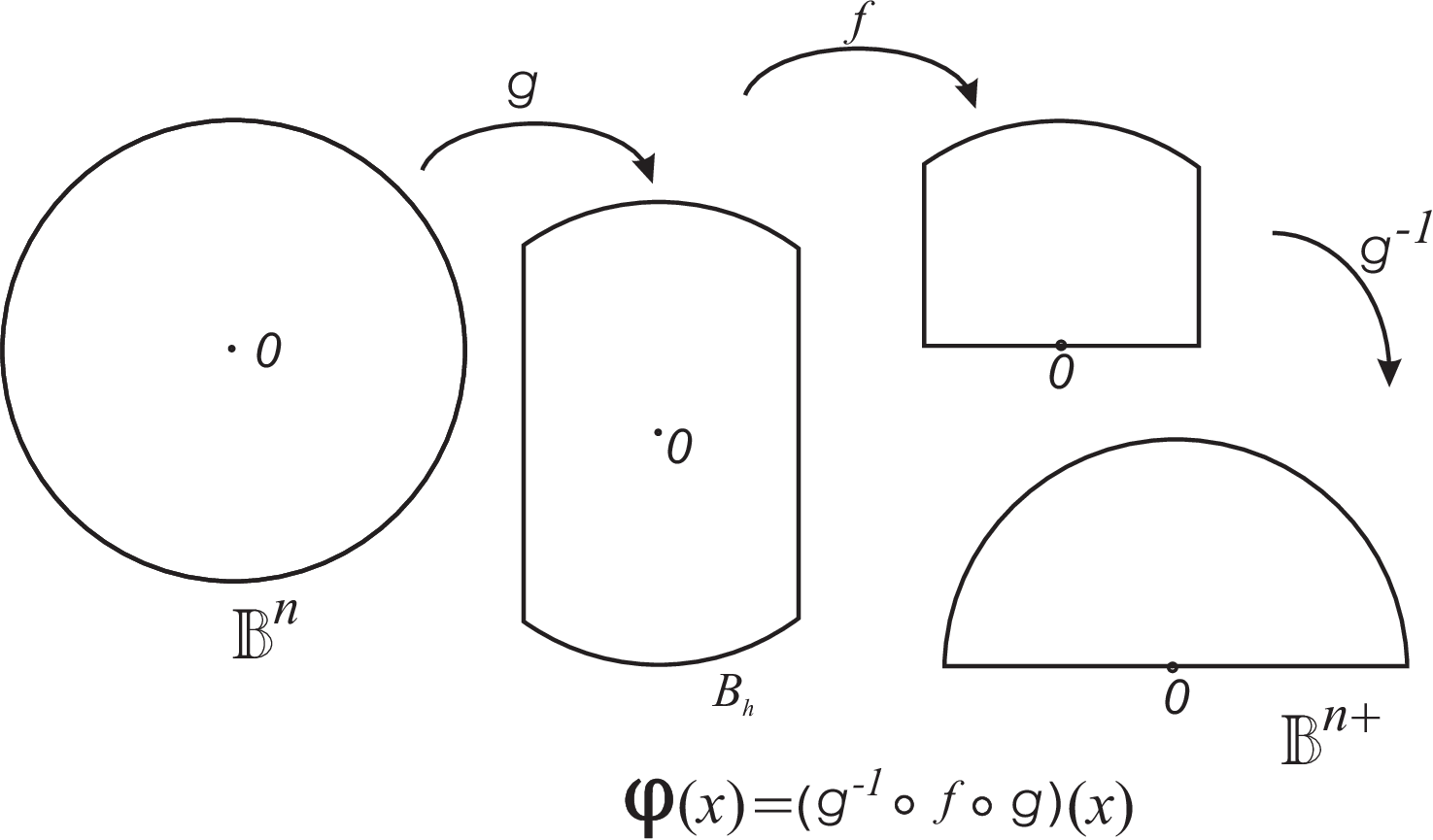}
  \caption{Illustration for Example~\ref{ex3}.}\label{fig3}
 \end{figure}

\medskip
Set $F(x):=(f^{\,-1}_0\circ\varphi)(x)$ and
$F_m(x):=\frac{m}{m+1}(f^{\,-1}_0\circ\varphi)(x).$ Now, $F$ maps
$D={\Bbb B}^n$ onto $D_0$ and $F_m$ is a quasiconformal mapping of
${\Bbb B}^n$ onto $D_0\cap B(0, m/(m+1)).$ The fact that the
sequence of domains $D_m:=D_0\cap B(0, m/(m+1))$ is regular with
respect to some sequence of domains of cuts $d_k$ of the domains
$D_0$ may be verified in exactly the same way as in the previous
example.

\medskip
Obviously, $D_0\cap B(0, m/(m+1))$ converge to $D_0$ as to its
kernel. In addition, $F_m$ has a continuous extension to
$\overline{D}$ by the metric $\rho,$ because $F_m$ has a usual
(Euclidean) continuous extension to $\overline{D},$ which, in turn,
is ensured by the finite connectedness of the domain $D_0$ on its
boundary (see e.g. \cite[Theorem~4.2(3)]{Na$_1$}). In addition,
$D_m$ are compact subdomains of $D_0$ and metric $\rho$ in
$\overline{D_0}_P$ is homeomorphic to the Euclidean metric onto
compact sets of $D_0.$

For $y_0:=(-1/4, 0,\ldots, 0, 0)$ we have that
$F^{\,-1}_m(y_0)\rightarrow F^{\,-1}(y_0)\in {\Bbb B}^n$ as
$m\rightarrow\infty,$ consequently, the relation $h(F_m^{\,-1}(E),
\partial D)\geqslant \delta$ holds for $E=\{y_0\}$ and
$\delta=h(F^{\,-1}(y_0), \partial {\Bbb B}^n)/2$ for sufficiently
large $m\in {\Bbb N}.$ In addition, $F_m$ are homeomorphisms and,
consequently, are open, discrete and closed. The
relations~(\ref{eq2*A})--(\ref{eqA2}) holds for $F_m$ with some
$Q(x)\equiv K=const$ because $F_m$ are quasiconformal with some
general coefficient $K$ (see e.g. \cite{Va}). Obviously, $Q\in
L^1(D_0)$ and, consequently, satisfies the conditions mentioned in
Theorem~\ref{th2}. Thus, all conclusions of Theorem~\ref{th2} hold.
By the same reasons, all of the conditions of Theorem~\ref{th1} hold
for some a continuum $E_2\subset D$ and $\delta_2>0$ and conclusions
of this theorem hold, as well.
\end{example}

\medskip
\begin{example}\label{ex4}
It is easy to find corresponding examples of mappings with branching
that satisfy the conditions of Theorems~\ref{th2} or~\ref{th1}. In
particular, in the notation of Example~\ref{ex2} we may put
$\widetilde{f}_m=f_m\circ \varphi$ and $\widetilde{f}=f\circ
\varphi,$ where $\varphi(z)=z^2.$ The domains $D, D_0, D_m, d_k$ are
not change in this case as well as a compactum $E$ mentioned in
Example~\ref{ex2}. The mapping $\varphi(z)=z^2$ is quasiregular and,
consequently, satisfies the condition~(\ref{eq2*A}) at any point
$y_0\in \overline{D_0}$ with $Q(y)=2$ (see \cite[Theorem~2.4.I,
Remark~2.5.I]{Ri}). Therefore, $\widetilde{f}_m$ also satisfy
(\ref{eq2*A})--(\ref{eqA2}) with the mentioned $Q.$ The mappings
$\widetilde{f}_m$ are open, discrete and closed and satisfy all the
conditions of Theorem~\ref{th2} (Theorem~\ref{th1}). Similarly, in
${\Bbb R}^n,$ $n\geqslant 3,$ we take $\varphi(x)=(z, r\cos
l\varphi, r\sin l\varphi),$ $\widetilde{f}_m=f_m\circ \varphi$ and
$\widetilde{f}=f\circ \varphi,$ where $l\in {\Bbb N},$ $l\geqslant
2.$ Observe that $K_O(x, \varphi)=l^{n-1}$ (see item~4 Ch.~I
in~\cite{Re}) and, correspondingly, $Q(x)=l^{n-1}$ in~(\ref{eq2*A})
(see \cite[Theorem~2.4.I, Remark~2.5.I]{Ri}). Therefore,
$\widetilde{f}_m$ also satisfy (\ref{eq2*A})--(\ref{eqA2}) with
$C\cdot Q,$ where $C$ is some constant and $Q=l^{n-1}.$ Obviously,
the mappings $\widetilde{f}_m$ are open, discrete and closed and
satisfy all the conditions of Theorem~\ref{th2} (Theorem~\ref{th1}).
\end{example}

\medskip
\begin{example}\label{ex5}
It is possible to construct corresponding families of mappings,
satisfying Theorem~\ref{th2} (Theorem~\ref{th1}) that have unbounded
characteristics. This can be done, for example, as follows. Let
$x_0\in D_0$ and $0<r_0<d(x_0, \partial D_0),$ where $D_0$ is a
domain from Example~\ref{ex2} for $n=2,$ and is a domain from
Example~\ref{ex3} for  $n\geqslant 3.$ Put
$h(x)=\frac{x}{|x|\log\frac{(r_0e)}{|x|}},$ $x\in B(x_0, r_0),$
$h(x_0)=x_0,$ $h|_{S(x_0, r_0)}=x.$ We denote $h_1(y):=h^{\,-1}(y).$
Then $h_1$ is defined in the ball $B(x_0, r_0)$ and $h_1(B(x_0,
r_0))=B(x_0, r_0).$ Reasoning similarly to
\cite[Proposition~6.3]{MRSY}, it may be shown that $h_1$ satisfies
the relations~(\ref{eq2*A})--(\ref{eqA2}) at any point $y_0\in
\overline{B(x_0, r_0)}$ for
$Q=Q(y)=\log^{n-1}\left(\frac{r_0e}{|y|}\right).$

Note that $Q\in L^1(B(x_0, r_0)).$ Indeed, by the Fubini theorem, we
calculate that
$$\int\limits_{B(x_0, r_0)}Q(y)\,dm(y)=\int\limits_0^{r_0}\int\limits_{S(0, r)}
\log^{n-1}\left(\frac{r_0e}{|y|}\right)d\mathcal{H}^{n-1}(y)dr=$$
$$=\omega_{n-1}\int\limits_0^1
r^{n-1}\log^{n-1}\left(\frac{r_0e}{r}\right)\,dr\leqslant
\omega_{n-1}(r_0e)^{n-1}\int\limits_0^1dr=\omega_{n-1}(r_0e)^{n-1}<\infty\,,$$
where $\omega_{n-1}$ denotes the area of the unit sphere ${\Bbb
S}^{n-1}$ in ${\Bbb R}^n.$

\medskip
Now, let for $n=2,$ in the notations of Examples~\ref{ex2}
and~\ref{ex3}
$$g_m(x)=\begin{cases}(h_1\circ f_m)(x)\,,& x\in f^{\,-1}_m(B(x_0, r_0))\,,\\
f_m(x)\,,& x\not\in f^{\,-1}_m(B(x_0, r_0))\end{cases}\,,$$
$$g(x)=\begin{cases}(h_1\circ f)(x)\,,& x\in f^{\,-1}(B(x_0, r_0))\,,\\
f(x)\,,& x\not\in f^{\,-1}(B(x_0, r_0))\end{cases}\,,$$
where $f_m$ and $f$ are from Example~\ref{ex2}, and for $n\geqslant
3$
$$G_m(x)=\begin{cases}(h_1\circ F_m)(x)\,,& x\in F^{\,-1}_m(B(x_0, r_0))\,,\\
f_m(x)\,,& x\not\in f^{\,-1}_m(B(x_0, r_0))\end{cases}\,,$$
$$G(x)=\begin{cases}(h_1\circ F)(x)\,,& x\in F^{\,-1}(B(x_0, r_0))\,,\\
F(x)\,,& x\not\in F^{\,-1}(B(x_0, r_0))\end{cases}\,,$$
where $F_m$ and $F$ are from Example~\ref{ex3}.
By the construction, $g_m$ and $G_m$ satisfy the
relations~(\ref{eq2*A})--(\ref{eqA2}) at any point $y_0\in
\overline{D_0}$ for $Q=C\cdot
Q(y)=\log^{n-1}\left(\frac{r_0e}{|y|}\right),$ where $C$ is some
constant. Observe that, the domains $f_m(D)$ and $F_m(D)$ are not
changed under applying to them the mapping $h_1,$ so that the
domains $g_m(D)$ and $G_m(D)$ are regular, as well. In addition,
since $h_1$ is a fixed mapping, the relations $h(g^{\,-1}_m(E),
\partial D)\geqslant \delta>0$ or $h(G^{\,-1}_m(E),
\partial D)\geqslant \delta>0$ also hold for some a compact $E\subset D_0,$
$\delta>0$ and all $m=1,2,\ldots .$ The mappings $g_m$ and $G_m,$
$m=1,2,\ldots ,$ are open, discrete and closed and satisfy all the
conditions of Theorem~\ref{th2} (Theorem~\ref{th1}).
\end{example}

\medskip
\begin{example}\label{ex6}
Let $\widetilde{D}$ be the unit square from which the sequence of
segments $I_k=\{z=(x, y)\in {\Bbb R}^2: x=1/k,\,\,0<y\leqslant
1/2\},$ $k=2,3,\ldots,$ is removed (see Figure~\ref{fig7}).
\begin{figure}
  \centering\includegraphics[scale=0.5]{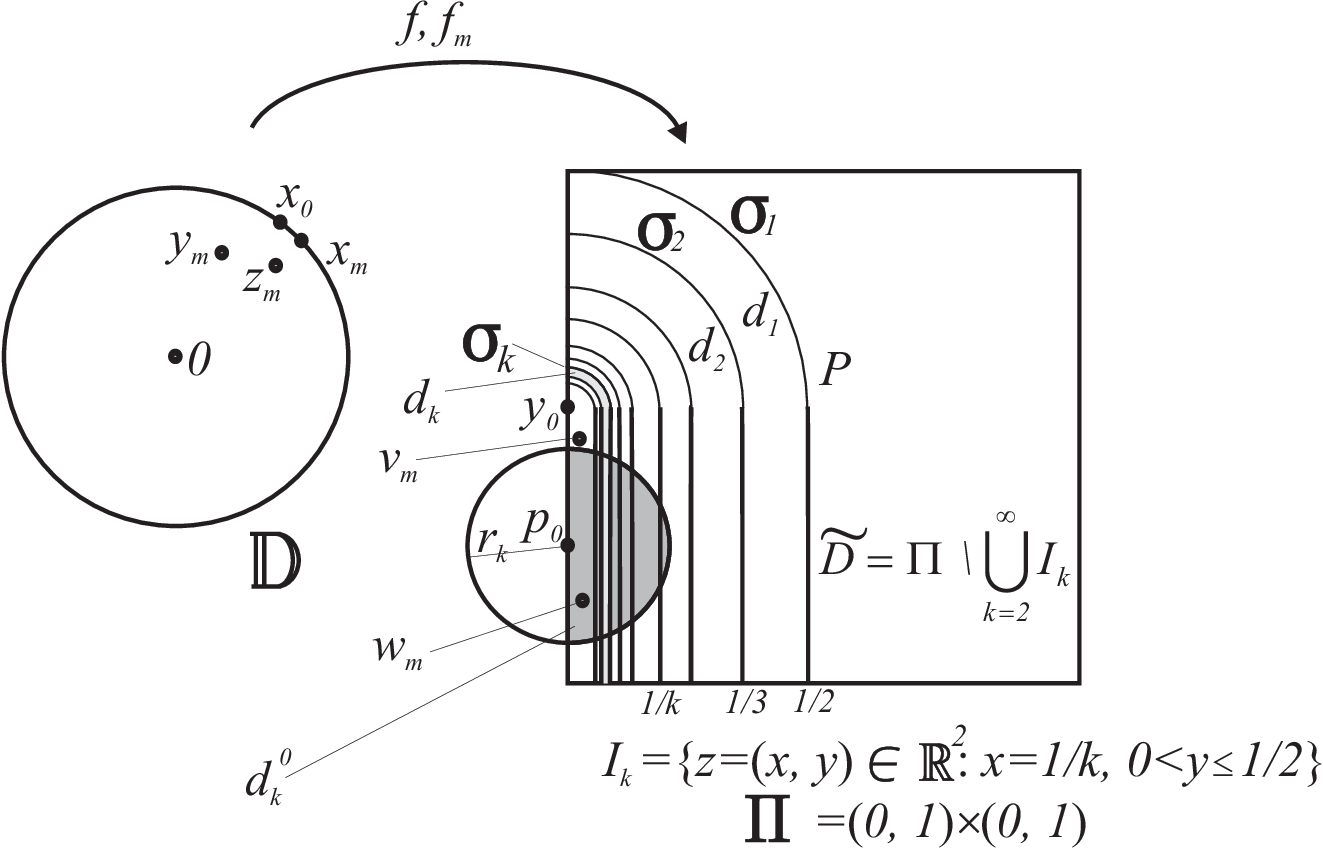}
  \caption{Illustration for Example~\ref{ex6}.}\label{fig7}
 \end{figure}
Consider the prime end $P$ in the domain $\widetilde{D},$ formed by
cuts
$$\sigma_m=\left\{z=y_0+\frac{e^{i\varphi}}{m+1},\,\, y_0=(0, 1/2),\,\, 0\leqslant
\varphi\leqslant \pi/2\right\}, \quad m=1,2,\ldots, .$$
It can be shown that the end $P$ is really prime. According to the
Riemannian mapping theorem, there exists a conformal mapping $f$ of
the unit disk ${\Bbb D}$ onto the domain $\widetilde{D}$ and by the
Caratheodory theorem, the prime end $P$ corresponds to some
point~$x_0\in\partial {\Bbb D}$ so that $C(f, x_0)=I(P),$
see~\cite[Theorem~9.4]{CL}. By the same theorem, the above
correspondence is one-to-one.

\medskip
Let $D_m$ be a domain obtained from $\widetilde{D}$ by the removing
of the segments $I_1, I_2,\ldots, I_m, I_{m+1},$ $m=1,2,\ldots .$
Again by the Riemannian theorem, there is a mapping $f_m$ of ${\Bbb
D}$ onto $D_m,$ $m=1,2,\ldots.$ Due to the additional fractional
linear transformation, we may consider that $f_m(0)=(3/4, 3/4).$ By
the reasons mentioned above, $f^{\,-1}_m$ has a continuous extension
$f_m:\overline{D_m}_P\rightarrow \overline{D}.$ Let $x_m\in \partial
{\Bbb D}$ be a point, $m=1,2,\ldots,$ such that $f^{\,-1}_m(P)=x_m.$
Due to this, one can chose sequences $v_m, w_m$ in $D_m$ such that
$v_m, w_m\rightarrow P$ as $m\rightarrow\infty,$ $|v_m-w_m|>1/4$ and
$f^{\,-1}_m(v_m)-x_m\rightarrow 0,$ $f^{\,-1}_m(w_m)-x_m\rightarrow
0$ as $m\rightarrow\infty.$ Due to the compactness of ${\Bbb
S}^1=\partial {\Bbb D},$ we may consider that $x_m\rightarrow x_0\in
\partial {\Bbb D}.$ Thus, the sequences $y_m:=f^{\,-1}_m(v_m)$ and
$z_m:=f^{\,-1}_m(w_m)$ converge to $x_0$ as $m\rightarrow\infty.$ On
the other hand, the sequence $D_m$ converges to $D_0:=[0, 1]\times
[0, 1]$ as its kernel. In addition, $f_m$
satisfy~(\ref{eq2*A})--(\ref{eqA2}) at any point $y_0\in
\overline{D_0}$ for $Q\equiv 1$ (see \cite[Theorem~8.1]{Va}). The
condition $h(f^{\,-1}_m((3/4, 3/4)), \partial {\Bbb
D})\geqslant\delta $ obviously holds for some $\delta>0$ and all
$m=1,2,\ldots$ because $f^{\,-1}_m((3/4, 3/4))=0$ by the
construction. Observe that, the sequence $\{f_m\}_{m=1}^{\infty}$ is
not equicontinuous as a sequence between metric spaces $({\Bbb D},
|\cdot|)$ and $(\overline{D_0}_P, \rho)$ because
$|f_m(x_m)-f_m(y_m)|=|v_m-w_m|\geqslant 1/4\not\rightarrow 0$ as
$m\rightarrow\infty;$ in addition, the pointwise convergence in
$\overline{D_0}$ is equivalent to the convergence by the metric
$\rho$ in prime ends space $(\overline{D_0}_P, \rho)$ because $D_0$
is a domain with a locally quasiconformal boundary, consequently, we
may set $|\rho(x)-\rho(y)|=|x-y|$ for $x,y\in \overline{D_0}_P.$
Thus, the conclusions of Theorem~\ref{th2} and~\ref{th1} do not
hold. The reasons of that are the following: mappings $f_m$ have no
continuous boundary extension $f_m:\overline{D}\rightarrow
\overline{D_0}_P$ by the metric $\rho$ in $\overline{D_0}_P,$ in
addition, any sequence of domains $d^0_k,$ $k=1,2,\ldots ,$ (in
particular, the sequence $d^0_k$ formed by sufficiently small
circles centered at some a point $p_0\in E_{D_0}$ of the radius
$r_k>0,$ $r_k\rightarrow 0$ as $k\rightarrow\infty$ of the axes
$O_y$) formes non-connected intersection with $D_m$ for all
$m=1,2,\ldots $ and infinitely many $k=1,2,\ldots .$ Any such a
sequence of domains $d^0_k,$ $k=1,2,\ldots,$ corresponds one and
only one prime end in $E_{D_0},$ see
e.g.~\cite[Theorem~4.1]{Na$_2$}.
\end{example}

\medskip\medskip\medskip
{\bf Funding.} The paper was supported by the National Research
Foundation of Ukraine, Project number 2025.07/0014, Project name:
``Modern problems of Mathematical Analysis and Geometric Function
Theory''.

\medskip
\noindent{{\bf Zarina Kovba} \\
Zhytomyr Ivan Franko State University,  \\
40 Velyka Berdychivs'ka Str., 10 008  Zhytomyr, UKRAINE \\
mazhydova@gmail.com

\medskip
{\bf \noindent Evgeny Sevost'yanov} \\
{\bf 1.} Zhytomyr Ivan Franko State University,  \\
40 Velyka Berdychivs'ka Str., 10 008  Zhytomyr, UKRAINE \\
{\bf 2.} Institute of Applied Mathematics and Mechanics\\
of NAS of Ukraine, \\
19 Henerala Batyuka Str., 84 116 Slov'yans'k,  UKRAINE\\
esevostyanov2009@gmail.com

\end{document}